\documentclass[11pt]{article}
\usepackage{amscd}
\usepackage{amsfonts}
\usepackage{amsmath}
\usepackage{amssymb}
\usepackage{amsthm}
\usepackage{bbm}
\usepackage{CJK}
\usepackage{fancyhdr}
\usepackage{graphicx}
\usepackage{hyperref}
\usepackage{indentfirst}
\usepackage{latexsym}
\usepackage{mathrsfs}
\usepackage{xypic}

\usepackage[top=1in,bottom=1in,left=1.25in,right=1.25in]{geometry}
\def\<{\langle}
\def\>{\rangle}
\def\a{\alpha}
\def\b{\beta}

\def\c{\cdot}

\def\D{\Delta}

\def\g{\gamma}

\def\lr{\longrightarrow}

\def\o{\otimes}

\def\r{\rho}

\def\v{\varepsilon}

\date{}
\begin{document}
\renewcommand{\baselinestretch}{1.2}
\renewcommand{\arraystretch}{1.0}
\title{{\bf  Braided monoidal categories and Doi Hopf modules for monoidal Hom-Hopf algebras   }
\date{}}
\author{{\bf Shuangjian Guo$^{1}$, Xiaohui Zhang$^{2}$
        and Shengxiang Wang$^{3}$\footnote
        {Correspondence: S X Wang: wangsx-math@163.com.} }\\
1.~ School of Mathematics and Statistics, Guizhou University of Finance and\\ Economics, Guizhou  550025, China\\
2.~Department of Mathematics, Southeast  University,
Nanjing 210096, China\\
 3.~ School of Mathematics and Finance, Chuzhou University,
 Chuzhou 239000,  China}
 \maketitle

 {\bf Abstract} We first introduce the notion of Doi Hom-Hopf modules
  and find the sufficient condition for the category of  Doi Hom-Hopf modules to be monoidal. 
Also we obtain the condition for the monoidal Hom-algebra  and monoidal Hom-coalgebra to be  monoidal Hom-bialgebras. 
Second, we  give the maps between the underlying monoidal Hom-Hopf algebras,
  Hom-comodule algebras and Hom-module coalgebras give rise to functors between
   the category of  Doi Hom-Hopf modules and study tensor identities for monodial categories of Doi Hom-Hopf modules.
   Furthermore, we construct a braiding on the category of  Doi Hom-Hopf modules.
  Finally, as an application of our theory, we consider the braiding on the category of 
   Hom-modules, the category of  Hom-comodules and the category of Hom-Yetter-Drinfeld modules respectively.
\medskip

\medskip

\textbf {MR(2010)Subject Classification} 16T05

\section {Introduction}
\def\theequation{1. \arabic{equation}}
\setcounter{equation} {0} \hskip\parindent

The category $_{ A}\mathcal{M}(H)^{C} $ of Doi-Hopf modules  was introduced
in \cite{D83}, where $H$ is a Hopf algebra, $A$ a right  $H$-comodule algebra and $C$ a left $H$-module coalgebra. It is the category of the modules over the algebra
$A$  which are also comodules over the coalgebra $C$ and satisfy certain compatibility
condition involving $H$. The study of $_{ A}\mathcal{M}(H)^{C} $ turned out to be very useful: it was
shown in \cite{D83} that many categories  such as the
module and comodule categories over bialgebras, the Hopf modules category \cite{S69},
and the Yetter-Drinfeld category \cite{RT93} are special cases of
$_{ A}\mathcal{M}(H)^{C} $. For a further study of Doi-Hopf modules, we refer to \cite{CR95}, \cite{CVZ98}.
In \cite{CM97}, they proved that  Yetter-Drinfel'd modules are special
cases of Doi-Hopf modules,  therefore the category of Yetter-Drinfel'd modules
is  a Grothendieck category.

Hom-algebras and Hom-coalgebras were introduced by Makhlouf and Silvestrov in
\cite{MS08} as generalizations of ordinary algebras and coalgebras in the following sense: the
associativity of the multiplication is replaced by the Hom-associativity and similar for
Hom-coassociativity. They also defined the structures of Hom-bialgebras and Hom-Hopf
algebras, and described some of their properties extending properties of ordinary bialgebras
and Hopf algebras in \cite{MS09} and \cite{MS10}.
Recently, many more properties and structures of Hom-Hopf algebras have been developed,
see \cite{CZ13},\cite{CWZ13}, \cite{CZ14}, \cite{GC14}, \cite{LS14},\cite{WG},\cite{WCZ12}  and references cited therein.

Caenepeel and Goyvaerts \cite{CG11} studied Hom-bialgebras and Hom-Hopf algebras from
a categorical view point, and called them monoidal Hom-bialgebras and monoidal Hom-
Hopf algebras respectively, which are slightly different from the above Hom-bialgebras
and Hom-Hopf algebras.
In \cite{MP14}, Makhlouf and Panaite defined Yetter-Drinfeld modules over Hom-bialgebras
 and shown that Yetter-Drinfeld modules over a Hom-bialgebra with bijective structure map provide solutions of the Hom-Yang-Baxter equation.
Also Liu and Shen \cite{LS14} studied Yetter-Drinfeld modules over monoidal Hom-bialgebras and called them Hom-Yetter-Drinfeld modules,
 and shown that the category of Hom-Yetter-Drinfeld modules is a braided monoidal categories.
Chen and Zhang  \cite{CZ14} defined the category of Hom-Yetter-Drinfeld modules in a slightly different way to \cite{LS14},
 and shown that it is a full monoidal subcategory of the left center of left Hom-module category.

In this paper, we will discuss the following questions: how do we define a Doi Hopf module for monoidal Hom-Hopf algebras (we will call it a Doi Hom-Hopf module)
 such that the category of  Doi Hom-Hopf modules is monoidal?
 In section 3, we will show that it is sufficient that $(A, \b)$ and $(C,\g)$ are monoidal Hom-bialgebras with some extra conditons. 
 As an example, we consider the category of Hom-Yetter-Drinfeld modules, which is well known to be a monoidal category from \cite{LS14},
  this category is a special of our theory. 
  
  In section 4, we  give the maps between the underlying monoidal Hom-Hopf algebras, Hom-comodule algebras and Hom-module coalgebras give rise to functors between the category of  Doi Hom-Hopf modules. If the maps between the monodial Hom-algebras and monodial Hom-coalgebras are both monoidal Hom-bialgerbas maps, then we obtain functors monodial categories, and study tensor identities for monodial categories of Doi Hom-Hopf modules.
   As an application, we prove that the category of  Doi Hom-Hopf modules has enough injective objects.

Suppose that we have a monoidal  category of  Doi Hom-Hopf modules. 
How do we define a braiding on this category? In section 5, we point out this comes down to giving a twisted convolution inverse map $\mathscr{Q}: \,\, C\o C\rightarrow A\o A$ satisfying the complicated compatibility conditions, as an application we consider the braiding on the category of  Hom-modules, the category of  Hom-comodules and the category of Hom-Yetter-Drinfeld modules respectively.

\section{Preliminaries}
Throughout this paper we work over a commutative ring $k$, we recall from \cite{CG11} and \cite{GC14}
for some informations about Hom-structures which are needed in what follows.

Let $\mathcal{C}$ be a category. We introduce a new category
$\widetilde{\mathscr{H}}(\mathcal{C})$ as follows: objects are
couples $(M, \mu)$, with $M \in \mathcal{C}$ and $\mu \in
Aut_{\mathcal{C}}(M)$. A morphism $f: (M, \mu)\rightarrow (N, \nu )$
is a morphism $f : M\rightarrow N$ in $\mathcal{C}$ such that $\nu
\circ f= f \circ \mu$.

Let $\mathscr{M}_k$ denotes the category of $k$-modules.
~$\mathscr{H}(\mathscr{M}_k)$ will be called the Hom-category
associated to $\mathscr{M}_k$. If $(M,\mu) \in \mathscr{M}_k$, then
$\mu: M\rightarrow M$ is obviously a morphism in
~$\mathscr{H}(\mathscr{M}_k)$. It is easy to show that
~$\widetilde{\mathscr{H}}(\mathscr{M}_k)$ =
(~$\mathscr{H}(\mathscr{M}_k),~\otimes,~(I, I),~\widetilde{a},
~\widetilde{l},~\widetilde{r}))$ is a monoidal category by
Proposition 1.1 in \cite{CG11}: the tensor product of $(M,\mu)$ and $(N,
\nu)$ in ~$\widetilde{\mathscr{H}}(\mathscr{M}_k)$ is given by the
formula $(M, \mu)\otimes (N, \nu) = (M\otimes N, \mu \otimes \nu)$.

Assume that $(M, \mu), (N, \nu), (P,\pi)\in
\widetilde{\mathscr{H}}(\mathscr{M}_k)$. The associativity and unit
constraints are given by the formulas
\begin{eqnarray*}
\widetilde{a}_{M,N,P}((m\o n)\o p)=\mu(m)\o (n\o \pi^{-1}(p)),\\
\widetilde{l}_{M}(x\o m)=\widetilde{r}_{M}(m\o x)=x\mu(m).
\end{eqnarray*}
An algebra in $\widetilde{\mathscr{H}}(\mathscr{M}_k)$ will be called a monoidal Hom-algebra.

{\bf Definition 2.1.} A monoidal Hom-algebra  is an
object $(A, \alpha)\in
\widetilde{\mathscr{H}}(\mathscr{M}_k)$ together with a $k$-linear
map $m_A: A\o A\rightarrow A$ and an element $1_A\in A$ such that
\begin{eqnarray*}
\alpha(ab)=\alpha(a)\alpha(b); ~\a(1_A)=1_A,\\
\alpha(a)(bc)=(ab)\alpha(c); ~a1_A=1_Aa=\a(a),
\end{eqnarray*}
for all $a,b,c\in A$. Here we use the notation $m_A(a\o b)=ab$.

{\bf Definition 2.2.} A monoidal Hom-coalgebra  is
an object $(C,\g)\in \widetilde{\mathscr{H}}(\mathscr{M}_k)$
together with $k$-linear maps $\Delta:C\rightarrow C\otimes C,~~~
\D(c)=c_{(1)}\o c_{(2)}$ (summation implicitly understood) and
$\varepsilon:C\rightarrow k$ such that
\begin{eqnarray*}
\D(\g(c))=\g(c_{(1)})\o \g(c_{(2)}); ~~~\varepsilon(\g(c))=\varepsilon(c),
\end{eqnarray*}
and
\begin{eqnarray*}
\g^{-1}(c_{(1)})\otimes c_{(2)(1)}\otimes c_{(2)(2)}=c_{(1)(1)}\otimes c_{(1)(2)}\otimes\g^{-1}(c_{(2)}),~\varepsilon(c_{(1)})c_{(2)}=\varepsilon(c_{(2)})c_{(1)}=\g^{-1}(c)
\end{eqnarray*}
for all $c\in C$.

{\bf Definition 2.3.} A monoidal Hom-bialgebra
$H=(H,\a,m,\eta, \Delta,\varepsilon)$  is a bialgebra in
the symmetric monoidal category
$\widetilde{\mathscr{H}}(\mathscr{M}_k)$. This means that $(H, \a,
m,\eta)$ is a monoidal Hom-algebra, $(H,\alpha,\Delta,\varepsilon)$ is a monoidal Hom-coalgebra
and that $\D$ and $\v$ are morphisms of Hom-algebras, that is,
\begin{eqnarray*}
\Delta(ab)=a_{(1)}b_{(1)}\otimes a_{(2)}b_{(2)}; ~~~\Delta(1_H)=1_H\otimes 1_H,\\
\varepsilon(ab)=\varepsilon(a)\varepsilon(b),~\varepsilon(1_H)=1_H.
\end{eqnarray*}

{\bf Definition 2.4.} A monoidal Hom-Hopf algebra is
a monoidal Hom-bialgebra $(H, \alpha)$ together with a linear map
$S:H\rightarrow H$ in $\widetilde{\mathscr{H}}(\mathscr{M}_k)$ such
that
$$S\ast I=I\ast S=\eta\varepsilon,~S\alpha=\alpha S.$$

{\bf Definition 2.5.} Let $(A,\alpha)$ be a monoidal
Hom-algebra.  A right $(A,\alpha)$-Hom-module is an object
$(M,\mu)\in \widetilde{\mathscr{H}}(\mathscr{M}_k) $ consists of a
$k$-module and a linear map $\mu:M\rightarrow M$ together with a
morphism $\psi:M\otimes A\rightarrow M, \psi(m\c a)=m\c a$, in
$\widetilde{\mathscr{H}}(\mathscr{M}_k) $  such that
\begin{eqnarray*}
(m\c a)\c\alpha(b)=\mu(m)\c (ab);~~~~ m\c 1_A=\mu(m),
\end{eqnarray*}
for all $a\in A$ and $m\in M$. The fact that $\psi\in \widetilde{\mathscr{H}}(\mathscr{M}_k)$ means that
\begin{eqnarray*}
\mu(m\c a)=\mu(m)\c\alpha(a).
\end{eqnarray*}
A morphism $f: (M, \mu)\rightarrow (N, \nu)$ in $\widetilde{\mathscr{H}}(\mathscr{M}_k)$ is called right $A$-linear if it preserves the $A$-action, that is, $f(m\c a)=f(m)\c a$. $\widetilde{\mathscr{H}}(\mathscr{M}_k)_A$ will denote the category of right $(A, \a)$-Hom-modules and $A$-linear morphisms.

{\bf Definition 2.6.} Let $(C,\g)$ be a monoidal
Hom-coalgebra. A right $(C,\g)$-Hom-comodule  is an object
$(M,\mu)\in \widetilde{\mathscr{H}}(\mathscr{M}_k)$  together with a
$k$-linear map $\rho_M: M\rightarrow M\otimes C$ notation
$\rho_M(m)=m_{[0]}\o m_{[1]}$ in
$\widetilde{\mathscr{H}}(\mathscr{M}_k)$ such that
\begin{eqnarray*}
m_{[0][0]}\otimes (m_{[0][1]}\otimes
\g^{-1}(m_{[1]}))=\mu^{-1}(m_{[0]})\otimes \D_C(m_{[1]});
~m_{[0]}\varepsilon(m_{[1]})=\mu^{-1}(m),
\end{eqnarray*}
for all $m\in M$.  The fact that $\rho_M\in \widetilde{\mathscr{H}}(\mathscr{M}_k)$ means that
\begin{eqnarray*}
\rho_M(\mu(m))=\mu(m_{[0]})\otimes \g(m_{[1]}).
\end{eqnarray*}
Morphisms of right $(C, \g)$-Hom-comodule are defined in the obvious way. The
category of right $(C, \g)$-Hom-comodules will be denoted by $\widetilde{\mathscr{H}}(\mathscr{M}_k)^C$.

{\bf Definition 2.7.}  Let $(H,\alpha)$  be a monoidal Hom-Hopf algebra.  A monoidal Hom-algebra  $(A,\beta)$ is called  a right
$(H,\alpha)$-Hom-comodule algebra if $(A,\beta)$ is a right $(H,\alpha)$ Hom-comodule with coaction $\rho_A: A\rightarrow A\o H,
 ~~\rho_A(a)= a_{[0]}\o a_{[1]}$
 such that the following conditions satisfy:
 \begin{eqnarray*}
 \rho_{A} ( ab)= a_{[0]} b_{[0]} \otimes a_{[1]}b_{[1]},
 \rho_{A} ( 1_A)=1_A\o 1_H,
 \end{eqnarray*}
for all $a,b \in A$.

\section {Making the category of Doi Hom-Hopf modules into a monoidal category}
\def\theequation{3. \arabic{equation}}
\setcounter{equation} {0} \hskip\parindent

{\bf Definition 3.1.} Let $(H,\alpha)$ be a monoidal Hom-Hopf algebra. A monoidal Hom-coalgebra $(C, \g)$ is called a left $(H,\alpha)$-Hom-module coalgebra, if $(C,\g)$ is a left $(H,\alpha)$-Hom-module with action $\phi: H \otimes C\rightarrow C$,
$\phi (h \otimes c) = h \cdot c$ such that the following conditions hold:
\begin{eqnarray*}
&&\Delta_C(h \cdot c)=h _{(1)} \cdot  c _{(1)} \otimes h _{(2)} \cdot  c _{(2)},\\
&&\varepsilon_C(h \cdot c) = \varepsilon_C(c)\varepsilon_H(h),
\end{eqnarray*}
for all  $c \in C$ and $g, h \in H$.

A  Doi Hom-Hopf datum is a triple ($H, A, C$), where $H$ is a monoidal Hom-Hopf algebra, $A$
a  right $(H, \a)$-Hom comodule algebra and $(C, \g)$ a  left  $(H, \a)$-Hom module coalgebra.

{\bf Definition 3.2.}  Given a  Doi Hom-Hopf datum ($H, A, C$).  A Doi Hom-Hopf module $(M,\mu)$
is a left $(A,\b)$-Hom-module which is also a right
$(C,\g)$-Hom-comodule with the  coaction structure $ \rho _{M} : M
\rightarrow M \otimes  C $ defined by $\rho _{M}(m)=m_{[0]}\o
m_{[1]}$ such that the following compatible condition holds: for all
$m\in M $ and $a\in A$,
\begin{eqnarray*}
&&\rho _{M}( a\c m ) = a _{[0]} \c m _{[0]} \otimes  a _{[1]}\c m _{[1]}.
\end{eqnarray*}

A morphism between left-right Doi Hom-Hopf modules is a $k$-linear map
which is a morphism in the categories ${}_A\widetilde{\mathscr{H}}(\mathscr{M}_k)$ and $\widetilde{\mathscr{C}}(\mathscr{M}_k)^C$ at the same time.  $ _A{}\widetilde{\mathscr{H}}(\mathscr{M}_k)(H)^{C}$ will denote  the category of left-right Doi Hom-Hopf modules and morphisms
between them.

Now suppose that $C$ and $A$ are both monoidal Hom-bialgebras.

{\bf Proposition 3.3.}
 Let $(M,\mu) \in {}_A{}\widetilde{\mathscr{H}}(\mathscr{M}_k)(H)^{C}$,
 $(N,\nu)\in {}_A{}\widetilde{\mathscr{H}}(\mathscr{M}_k)(H)^{C}$. Then we have
 $M \o N \in {}_A{}\widetilde{\mathscr{H}}(\mathscr{M}_k)(H)^{C}$ with structures:
\begin{eqnarray*}
&& a\c (m \o n)=a_{(1)}\c m \o a_{(2)}\c n,\\
&& \r _{M\o N}(m\o n )=m_{[0]}\o
 n_{[0]}\o m_{[1]}n_{[1]}
\end{eqnarray*}
if and only if the following condition holds:
\begin{eqnarray}
 a_{(1)[0]} \o a_{(2)[0]}\o (a_{(1)[1]}\c c)(a_{(2)[1]}\c d) =a_{[0](1)}\o a_{[0](2)}
  \o a_{[1]}\circ (cd),
\end{eqnarray}
for all $a\in A$ and $c, d\in C$.
Furthermore, the category $\mathcal{C}={}_A{}\widetilde{\mathscr{H}}(\mathscr{M}_k)(H)^{C}$ is a monoidal category.

{\bf Proof. } It is easy to see that $M\o N$ is  a left $(A, \b)$-module
 and that $M\o N$ is  a right
 $(C, \g)$-comodule. Now we check that the compatibility
  condition holds:
\begin{eqnarray*}
&&\r _{M\o N}(a\c (m\o n))\\
&=&(a_{(1)}\c m)_{[0]} \o (a_{(2)}\c n)_{[0]}\o
(a_{(1)}\c m)_{[1]} (a_{(2)}\c n)_{[1]} \\
&=& a_{(1)[0]}\c m_{[0]} \o (a_{(2)[0]}\c n_{[0]})\o (a_{(1)[1]}\c m_{[1]})(a_{(2)[1]}\c n_{[1]})\\
&\stackrel{(3.1)}{=}& a_{[0](1)}\c m_{[0]} \o (a_{[0](2)}\c n_{[0]})\o a_{[1]}\c (m_{[1]}n_{[1]})\\
&=&a_{[0]}\c ( m_{[0]}\o  n_{[0]})\o a_{[1]}\c (m_{[1]}n_{[1]}).
\end{eqnarray*}
So  $M \o N \in {}_A{}\widetilde{\mathscr{H}}(\mathscr{M}_k)(H)^{C}$.

 Conversely,  one can easily check that  $A\o C\in
 {}_A{}\widetilde{\mathscr{H}}(\mathscr{M}_k)(H)^{C}$,
 let $m=1\o c$ and $n=1\o d$ for any $c, d\in C$ and easily
 get Eq.(3.2).

Furthermore, $k$ is an object in ${}_A{}\widetilde{\mathscr{H}}(\mathscr{M}_k)(H)^{C}$ with
structures:
 $$
a\c x=\v _A(a)x, \quad\r (x)=x\o 1_C,
 $$
for all $x\in k$ if and only if the following condition holds:
\begin{equation}
\v _A(a)1_C=\v _A(a_{(0)})(a_{(1)}\c 1_C),
\end{equation}
for all $a\in A$.
Then it is easy to get that $(\mathcal{C}={}_A{}\widetilde{\mathscr{H}}(\mathscr{M}_k)(H)^{C}, \o, k, a, l, r)$ is a monoidal category,
where $a,l,r$ are given by the formulas:
\begin{eqnarray*}
&&\widetilde{a}_{M,N,P}((m\o n)\o p)=\mu(m)\o (n\o \pi^{-1}(p)),\\
&&\widetilde{l}_{M}(x\o m)=\widetilde{r}_{M}(m\o x)=x\mu(m),
\end{eqnarray*}
for $(M, \mu), (N, \nu), (P,\pi)\in \mathcal{C}$.
$\hfill \Box$

We call $G=(H,A,C)$ a \emph{monoidal Doi Hom-Hopf datum} if $G$ is a Doi Hom-Hopf datum, and $A, C$  are Hom-bialgebras with the additional compatibility relations Eq.(3.1) and Eq.(3.2).

We will give an example for the monodial category ${}_A{}\widetilde{\mathscr{H}}(\mathscr{M}_k)(H)^{C}$. First, we give the definition of Yetter-Drinfeld modules over a monoidal Hom-Hopf algebra,
which is also introduced by Liu and Shen in \cite{LS14} similarly.

{\bf Definition 3.4.} Let $(H, \a)$ be a monoidal Hom-Hopf algebra. A left-right $(H, \a)$-Hom-Yetter-Drinfeld module is an object $(M, \mu)$ in $\widetilde{\mathscr{H}}(\mathscr{M}_k)$, such that $(M, \mu)$ a left $(H, \a)$-Hom-module  and a  right $(H, \a)$-Hom-comodule  with the following compatibility condition:
\begin{eqnarray}
h_{(1)}\c m_{[0]}\o h_{(2)} m_{[1]}=\mu((h_{(2)}\c \mu^{-1}(m))_{[0]})\o ( h_{(2)}\c \mu^{-1}(m))_{[1]}h_{(1)}
\end{eqnarray}
for all $h\in H$ and $m\in M$. We denote by $_H{}\mathscr{HYD}^H$ the category of left-right $(H, \a)$-Hom-Yetter-Drinfeld
modules, morphisms being left $(H, \a)$-linear right $(H, \a)$-colinear maps.

{\bf Proposition 3.5.}  One has that Eq. (3.3) is equivalent to the following equation:
\begin{eqnarray}
\rho(h\c m)= \a(h_{(2)(1)})\c m_{[0]}\o (h_{(2)(2)}\a^{-1}(m_{[1]}))S^{-1}(h_{(1)}),
\end{eqnarray}
for all $h\in H$ and $m\in M$.

{\bf Proof.} For one thing, we compute
\begin{eqnarray*}
&&\mu((h_{(2)} \c \mu^{-1}(m))_{[0]}) \o ((h_{(2)} \c \mu^{-1}(m))_{[1]})h_{(1)} \\
&\stackrel {(3.4)}{=}& \mu(\a(h_{(2)(2)(1)}) \c \mu^{-1}(m_{[0]})) \o ((h_{(2)(2)(2)}\a^{-2}(m_{[1]}))S^{-1}(h_{(2)(1)}))h_{(1)}\\
&=&\a(h_{(2)(1)}) \c m_{[0]} \o (h_{(2)(2)}\a^{-1}(m_{[1]}))(S^{-1}(h_{(1)(2)})h_{(1)(1)})\\
&=&h_{(1)}\c m_{[0]} \o h_{(2)} m_{[1]}.
\end{eqnarray*}

For another, we have
\begin{eqnarray*}
&&h_{(2)(1)}\c m_{[0]} \o (h_{(2)(2)} m_{[1]})S^{-1}(h_{(1)}) \\
&\stackrel {(3.3)}{=}& \mu((h_{(2)(2)}\c \mu^{-1}(m))_{[0]})\o (( h_{(2)(2)}\c \mu^{-1}(m))_{[1]}h_{(2)(1)})S^{-1}(h_{(1)}) \\
&=& \mu((\a^{-1}(h_{(2)}) \c \mu^{-1}(m))_{[0]}) \o \a((\a^{-1}(h_{(2)}) \c \mu^{-1}(m))_{[1]})(h_{(1)(2)}S^{-1}(h_{(1)(1)}))\\
&=& \mu((\a^{-2}(h) \c \mu^{-1}(m))_{[0]}) \o \a^2((\a^{-2}(h) \c \mu^{-1}(m))_{[1]})\\
&=& (\a^{-1}(h) \c m)_{[0]} \o \a((\a^{-1}(h) \c m)_{[1]}),
\end{eqnarray*}
which implies Eq. (3.4). $\hfill \Box$

{\bf Theorem 3.6.} Let $(H, \a)$ be a monoidal Hom-Hopf algebra with a bijective antipode.

(1) $H$ can be made into a right $H^{op} \o H$-Hom-comodule algebra. The coaction $H\rightarrow H\o H^{op}\o H$ is given by the formula:
\begin{eqnarray*}
h\mapsto  \a(h_{(2)(1)})\o (S^{-1}(\a^{-1}(h_{(1)}))\o h_{(2)(2)}).
\end{eqnarray*}

(2) $H$ can be made into a left  $H^{op} \o H$-Hom module coalgebra. The action of $H^{op} \o H$ on $H$ is given by the formula:
\begin{eqnarray*}
(h\o k) \triangleright c =(k\a^{-1}(c))\a(h).
\end{eqnarray*}

(3) The category $_H\mathscr{HYD}^H$  of left-right Hom-Yetter-Drinfeld modules is isomorphic to a
category of Doi Hom-Hopf modules, namely $_H\widetilde{\mathscr{H}}(\mathscr{M}_k)(H^{op} \o H)^{H}$.

{\it Proof.} (1) We first prove that $H$ is a right $H^{op} \o H$-Hom comodule. For any $h\in H$,
\begin{eqnarray*}
&&(\a^{-1}\o \D_{H^{op} \o H})\rho_H(h)\\
&=&h_{(2)(1)}\o \D_{H^{op} \o H}( S^{-1}(\a^{-1}(h_{(1)}))\o h_{(2)(2)})\\
&=& h_{(2)(1)}\o   S^{-1}(\a^{-1}(h_{(1)(2)}))\o h_{(2)(2)(1)}\o  S^{-1}(\a^{-1}(h_{(1)(1)}))\o h_{(2)(2)(2)} \\
&=& \a(h_{(2)(1)(1)})\o   S^{-1}(\a^{-1}(h_{(1)(2)}))\o h_{(2)(1)(2)}\o  S^{-1}(\a^{-1}(h_{(1)(1)}))\o \a^{-1}(h_{(2)(2)}) \\
&=& \a^2(h_{(2)(2)(1)(1)})\o  S^{-1}(\a^{-1}(h_{(2)(1)}))\o \a(h_{(2)(2)(1)(2)})\o  S^{-1}(\a^{-2}(h_{(1)}))\o h_{(2)(2)(2)} \\
&=& \a^2(h_{(2)(1)(2)(1)})\o  S^{-1}(h_{(2)(1)(1)})\o \a(h_{(2)(1)(1)(2)})\o  S^{-1}(\a^{-2}(h_{(1)}))\o \a^{-1}(h_{(2)(2)})\\
&=&\rho(\a(h_{(2)(1)}))\o S^{-1}(\a^{-2}(h_{(1)}))\o \a^{-1}(h_{(2)(2)})\\
&=&(\rho_{H}\o \a^{-1})\rho_H(h)
\end{eqnarray*}
So $H$ is a right $H^{op} \o H$-Hom comodule.
We also have 
\begin{eqnarray*}
&&\rho(hg)= \a(h_{(2)(1)}g_{(2)(1)})\o (S^{-1}(h_{(1)}g_{(1)})\o \a^{-1}(h_{(2)(2)}g_{(2)(2)}))\\
&=&\a(h_{(2)(1)})\a(g_{(2)(1)})\o (S^{-1}(h_{(1)})S^{-1}(g_{(1)})\o \a^{-1}(h_{(2)(2)})\a^{-1}(g_{(2)(2)})) \\
&=&(\a(h_{(2)(1)})\o (S^{-1}(h_{(1)})\o \a^{-1}(h_{(2)(2)})))(\a(g_{(2)(1)})\o ( S^{-1}(g_{(1)})\o \a^{-1}(g_{(2)(2)}))))\\
&=& \rho_{H}(h)\rho_{H}(g).
\end{eqnarray*}
(2) Now  we  prove that $H$ is a $H^{op} \o H$-Hom comodule. For any $h, l, k, m, c\in H$, we have
\begin{eqnarray*}
&&(\a(l)\o \a(m))\triangleright  [(h\o k)\triangleright c]\\
&=&(\a(l)\o \a(m))\triangleright(k\a^{-1}(c))\a(h)\\
&=&[\a(m)[(\a^{-1}(k)\a^{-2}(c))h]]\a^2(l)\\
&=&[\a(m)[k(\a^{-2}(c))\a^{-1}(h))]]\a^2(l)\\
&=&\a(mk)[[\a^{-1}(c))h]\a(l)]\\
&=&\a(mk)[c(hl)]
=mk[c\a(hl)]\\
&=&(hl\o mk)\triangleright\a(c)=[(l\o m)(h\o k)]\triangleright\a(c),
\end{eqnarray*}
and this implies that $H$ is a $H^{op} \o H$-Hom module.

Using the fact that $(H, \a)$ is an $(H, \a)$-Hom-bimodule algebra, we can obtain
that $(H, \a)$ is a left $H^{op} \o H$-Hom-module coalgebra.

(3) Let  $(M, \c)$ be a left $(H, \a)$-module and $(M, \rho_M)$  a right
$(H, \a)$-comodule. Then $M\in {}_{H}\widetilde{\mathscr{H}}(\mathscr{M}_k)(H^{op} \o H)^{H}$ if and only if
\begin{eqnarray*}
\rho_M(h\c m)&=& \a(h_{(2)(1)})\c m_{[0]}\o (S^{-1}(\a^{-1}(h_{(1)}))\o h_{(2)(2)})\triangleright m_{[1]}\\
&=& \a(h_{(2)(1)})\c m_{[0]}\o (h_{(2)(2)}\a^{-1}(m_{[1]}))S^{-1}(h_{(1)}),
\end{eqnarray*}
for all $h\in H$ and $m\in M$. Thus $_{H}\widetilde{\mathscr{H}}(\mathscr{M}_k)(H^{op} \o H)^{H}$ is isomorphic to $_{H}\mathscr{HYD}^H$. $\hfill \Box$

{\bf Example 3.7.} Let $(H, \a)$ be a monoidal Hom-Hopf algebra, we have shown that the category of Doi Hom-Hopf modules $_H\widetilde{\mathscr{H}}(\mathscr{M}_k)(H^{op} \o H)^{H}$ and the category of Hom-Yetter-Drinfeld modules $_H\mathscr{HYD}^H$  are isomorphic. Recall from that \cite{LS14} that the category of Hom-Yetter-Drinfeld modules is a monoidal category; let us check that it is a special  case of Proposition 3.3.  Indeed,
take $A=H$ and $C=H^{op}$ as monoidal Hom-bialgebras. Let $a=h, c=k$ and $d=g$. Then the left-hand side amounts to
\begin{eqnarray*}
&&h_{[0](1)}\o h_{[0](2)}\o h_{[1]}\c(k\bullet g)\\
&=& \a(h_{(2)(1)(1)})\o \a(h_{(2)(1)(2)})\o (S^{-1}(\a^{-1}(h_{(1)}))\o h_{(2)(2)})\c(gk)\\
&=&\a(h_{(2)(1)(1)})\o \a(h_{(2)(1)(2)})\o [ h_{(2)(2)})\a^{-1}(gk)]S^{-1}(h_{(1)}).
\end{eqnarray*}
The right-hand side is
\begin{eqnarray*}
&&h_{(1)[0]}\o h_{(2)[0]}\o (h_{[1](1)}\c k)(h_{[1](2)}\c g)\\
&=& \a(h_{(1)(2)(1)})\o \a(h_{(2)(2)(1)})\o ((S^{-1}(\a^{-1}(h_{(1)(1)}))\o h_{(2)(2)(1)})\c k)\\
&&\hspace{7cm}\bullet ((S^{-1}(\a^{-1}(h_{(1)(2)}))\o h_{(2)(2)(2)})\c g)\\
&=& \a(h_{(1)(2)(1)})\o \a(h_{(2)(2)(1)})\o (( h_{(2)(2)(2)}\a^{-1}(g))S^{-1}(h_{(1)(2)}))\\
&&\hspace{7cm}(( h_{(2)(2)(1)}\a^{-1}(k))S^{-1}(h_{(1)(1)}))\\
&=& \a(h_{(1)(2)(1)})\o \a(h_{(2)(2)(1)})\o (( h_{(2)(2)(2)}\a^{-1}(g))[S^{-1}(\a^{-1}(h_{(1)(2)}))h_{(2)(2)(1)}])\\
&&\hspace{10cm}kS^{-1}(h_{(1)(1)})\\
&=& \a(h_{(1)(1)(2)})\o \a(h_{(2)(1)(2)})\o (( \a^{-1}(h_{(2)(2)})\a^{-1}(g))[S^{-1}(h_{(1)(1)(2)})\a^{-1}(h_{(2)(1)}]))\\
&&\hspace{10cm} kS^{-1}(\a(h_{(1)(1)(1)}))\\
&=& \a(h_{(2)(1)(1)})\o \a(h_{(2)(1)(2)})\o (( h_{(2)(2)(2)}\a^{-1}(g))[S^{-1}(h_{(2)(1)(1)})h_{(2)(1)(1)}])\\
&&\hspace{10cm} kS^{-1}(\a^{-1}(h_{(1)}))\\
&=& \a(h_{(2)(1)(1)})\o \a(h_{(2)(1)(2)})\o (( h_{(2)(2)}g) kS^{-1}(\a^{-1}(h_{(1)}))\\
&=& \a(h_{(2)(1)(1)})\o \a(h_{(2)(1)(2)})\o [( \a^{-1}(h_{(2)(2)})\a^{-1}(g)) k]S^{-1}(h_{(1)})\\
&=& \a(h_{(2)(1)(1)})\o \a(h_{(2)(1)(2)})\o [h_{(2)(2)} \a^{-1}(gk)]S^{-1}(h_{(1)}).
\end{eqnarray*}
$\hfill \Box$
\section{Tensor identities}
\def\theequation{4. \arabic{equation}}
\setcounter{equation} {0}

{\bf Theorem 4.1.} Given two  Doi Hom-Hopf datums $(H, A, C)$ and $(H', A', C')$, 
suppose that a morphism $\xi:(H,A,C)\rightarrow (H',A',C')$ consist three maps
$\varphi: H \rightarrow H'$, $\psi: A \rightarrow A '$ and $\phi: C \rightarrow C '$ which
are respectively monodial Hom-Hopf algebra, Hom-algebra and Hom-coalgebra maps satisfying
\begin{equation}
\phi (h \cdot c) = \varphi (h) \c \phi (c),
\end{equation}
\begin{equation}
\rho _{A ^{¡¯}} (\psi (a)) =\psi (a _{[0]}) \otimes \varphi (a _{[1]}),
\end{equation}
for all $c \in C$, $h \in H$ and $a \in A$. Then we have a functor
$F: {}_A\widetilde{\mathscr{H}}(\mathscr{M}_k)(H)^{C} \rightarrow {}_{A'}\widetilde{\mathscr{H}}(\mathscr{M}_k)(H')^{C'} $,
 defined as follows:
 $$
 F (M) = A' \otimes _{A} M,
 $$
 where $(A',\b')$ is a left $(A,\b)$-module via $\psi$ and with structure maps defined by
\begin{equation}
 b'\cdot( a ' \otimes _{A}m) = \b'^{-1}(b ')a ' \otimes _{A} \mu(m),
 \end{equation}
\begin{equation}
 \rho _{F(M)} (a' \otimes _{A} m) = a' _{[0]} \otimes _{A} m_{[0]} \otimes a' _{[1]}\c \phi( m _{[1]}),
\end{equation}
 for all $a', b' \in A'$ and  $m \in M$.

{\bf Proof.}
Let us show that $A' \otimes _{A} M$ is an object of $_{A'}\widetilde{\mathscr{H}}(\mathscr{M}_k)(H')^{C'} $. It is routine to check that $F (M)$ is a left $(A',\b')$-module. For this, we need to show that
$A' \otimes _{A} M$ is a right $(C',\g')$-comodule and satisfy the compatible condition,  for any $m \in M$ and $a', b' \in A'$, we have
\begin{eqnarray*}
\rho _{F(M)} (b'\c (a' \otimes _{A} m) )
&=& \rho _{F(M)}(\b'^{-1}(b ')a ' \otimes _{A} \mu(m))\\
&=& \b'^{-1}(b '_{[0]})a' _{[0]} \otimes _{A} \mu(m_{[0]}) \otimes [\b'^{-1}(b '_{[1]})a' _{[1]}]\c \phi( \g(m_{[1]}))\\
&=& b '_{[0]}[a' _{[0]} \otimes _{A} m_{[0]}] \otimes b '_{[1]}[a' _{[1]}\c \phi(m_{[1]})]\\
&=& b'\c( a' _{[0]} \otimes _{A} m_{[0]} \otimes a' _{[1]}\c \phi( m _{[1]}))=b'\rho _{F(M)}(a' \otimes _{A} m)
\end{eqnarray*}
i.e., the compatible condition holds. It remains  to prove that $A' \otimes _{A} M$ is a right $(C',\g')$-comodule.
 For any $m \in M$ and $a' \in A'$, we have
\begin{eqnarray*}
&&(\rho_{F(M)}\o id_{C'})\rho_{F(M)} (a' \otimes _{A} m)\\
&=& (\rho_{F(M)}\o id_C')( a' _{[0]} \otimes _{A} m_{[0]} \otimes a' _{[1]}\c \phi( m _{[1]}))\\
&=& [a' _{[0][0]} \otimes _{A}m_{[0][0]} \otimes a' _{[0][1]}\c \phi( m _{[0][1]})]\o a' _{[1]}\c \phi( m _{[1]})\\
&=& [\b'^{-1}(a' _{[0]}) \otimes _{A}\mu^{-1}(m_{[0]}) \otimes a' _{[1](1)}\c \phi( m _{[1](1)})]\o \a'(a' _{[1](2)})\c \phi( \g(m _{[1](2)}))\\
&=& a' _{[0]} \otimes _{A}m_{[0]} \otimes [a' _{[1](1)}\c \phi( m _{[1](1)})\o a' _{[1](2)}\c \phi( m _{[1](2)})]\\
&=& (id_{F(M)}\o \D_{C'})\rho_{F(M)}(a' \otimes _{A} m),
\end{eqnarray*}
and
\begin{eqnarray*}
&&(id_{F(M)}\o \varepsilon)\rho_{F(M)} (a' \otimes _{A} m)\\
&=&(id_{F(M)}\o \varepsilon) ( a' _{[0]} \otimes _{A} m_{[0]} \otimes a' _{[1]}\c \phi( m _{[1]}))\\
&=& a' _{[0]} \varepsilon(a' _{[1]})\otimes _{A} m_{[0]}\varepsilon(\phi( m _{[1]}))=a' \otimes _{A} m,
\end{eqnarray*}
as desired. And this completes the proof. $\hfill \Box$

{\bf Theorem 4.2.}
Under the assumptions of Theorem 4.1, we have a functor
$G:  {}_{A'}\widetilde{\mathscr{H}}(\mathscr{M}_k)(H')^{C'}\rightarrow {}_A\widetilde{\mathscr{H}}(\mathscr{M}_k)(H)^{C} $
which is right adjoint to $F$. $G$ is defined by
$$
G (M') = M' \Box _{C'} C,
$$
 with structure maps
\begin{equation}
a\c ( m' \otimes c )  = \beta (a _{[0]})\c  m' \otimes  a _{[1]}\c c,
\end{equation}
\begin{equation}
\rho _{G(M')} ( m'  \otimes c  ) =\mu'^{-1}( m') \otimes c _{(1)} \otimes \g(c _{(2)}) ,
\end{equation}
for all $a \in A $.

{\bf Proof.}
We first show that $G (M')$ is an object of $_A{}\widetilde{\mathscr{H}}(\mathscr{M}_k)(H)^{C} $. 
It is not hard to check that $G (M')$ is a left $(A,\b)$-module. Now we check that
$G (M')$ is a right $(C,\g)$-comodule and satisfy the compatible condition. For any $m' \in M'$ and $a\in A, c\in C$, we have
\begin{eqnarray*}
\rho _{G(M')} (a\c (m'\otimes c) )
&=& \rho _{G(M')}(\beta (a _{[0]})\c  m' \otimes  a _{[1]}\c c)\\
&=& a _{[0]}\c\mu'^{-1}( m') \otimes a _{[1](1)} \c c_{(1)}\otimes \a(a _{[1](2)}) \c\g(c _{(2)})\\
&=& \beta (a _{[0][0]})\c\mu'^{-1}( m') \otimes a _{[0][1]} \c c_{(1)}\otimes a _{[1]} \c\g(c _{(2)})\\
&=& a\c( \mu'^{-1}( m') \otimes c _{(1)} \otimes \g(c _{(2)}))=b'\rho _{G(M')}(m'\otimes c),
\end{eqnarray*}
i.e., the compatible condition holds. It remains  to prove that $M' \Box _{C'} C$ is a right $(C,\g)$-comodule. 
For any $m' \in M'$ and $a \in A$, we have
\begin{eqnarray*}
&&(\rho_{G(M')}\o id_{C'})\rho_{G(M')} (m' \otimes _{A} c)\\
&=& (\rho_{G(M')}\o id_{C'})(\mu'^{-1}( m') \otimes c _{(1)} \otimes \g(c _{(2)}))\\
&=& \mu'^{-2}( m') \otimes c _{(1)(1)} \otimes \g(c _{(1)(2)}) \otimes \g(c _{(2)})\\
&=&  \mu'^{-2}( m') \otimes \g^{-1}(c _{(1)}) \otimes \g(c _{(2)(1)}) \otimes \g^2(c _{(2)(2)})\\
&=&  \mu'^{-1}( m') \otimes c _{(1)} \otimes [\g(c _{(2)(1)}) \otimes \g(c _{(2)(2)})]\\
&=& (id_{G(M')}\o \D_{C})\rho_{G(M')}(m' \otimes  c),
\end{eqnarray*}
and
\begin{eqnarray*}
&&(id_{G(M')}\o \varepsilon)\rho_{G(M')} (m' \otimes  c)\\
&=&(id_{G(M')}\o \varepsilon) ( \mu'^{-1}( m') \otimes c _{(1)} \otimes \g(c _{(2)}))\\
&=& \mu'^{-1}( m') \otimes c _{(1)}\varepsilon(c _{(2)}) \otimes 1_C =m' \otimes  c,
\end{eqnarray*}
as required.

 $G(M')\in {}_A{}\widetilde{\mathscr{H}}(\mathscr{M}_k)(H)^{C} $ and the functorial properties can be checked in a straightforward way.
  Finally, we show that
 $G$ is a right adjoint to $F$. Take $(M, \mu) \in  {}_A{}\widetilde{\mathscr{H}}(\mathscr{M}_k)(H)^{C}$, define
 $\eta _{M}: M \rightarrow GF(M)=(M \otimes _{A} A')\square _{C'}C$ as follows: for all $m \in M$,
 $$
 \eta _{M}(m)=m_{[0]} \otimes _{A} 1_{A'} \otimes m _{[1]}.
 $$
It is easy to see  that $\eta _{M}\in   {}_A{}\widetilde{\mathscr{H}}(\mathscr{M}_k)(H)^{C}$.
Take $(M',\mu') \in  {}_{A'}\widetilde{\mathscr{H}}(\mathscr{M}_k)(H')^{C'}$, define $\delta _{M'}: FG(M')\rightarrow M'$,
where
$$
\delta _{M'}(m' \otimes c )\otimes _{A} a ')=\varepsilon _{C}(c) m'\cdot a ',
$$
It is easy to check that $\delta _{M'}$  is $(A,\b)$-linear and therefore $\delta _{M'}\in  {}_{A'}\widetilde{\mathscr{H}}(\mathscr{M}_k)(H')^{C'}$.
 We can also verify $\eta$ and $\delta$ defined above are all natural transformations and 
satisfy
$$
G(\delta _{M'})\circ \eta _{G(M')}=I,\    \ \delta _{F(M)}\circ F(\eta _{M})=I,
$$
for all $ M \in  {}_A{}\widetilde{\mathscr{H}}(\mathscr{M}_k)(H)^{C}$ and $M' \in  {}_{A'}\widetilde{\mathscr{H}}(\mathscr{M}_k)(H')^{C'}$.
And this completes the proof. $\hfill \Box$

 A morphism  $\xi=(\varphi, \psi, \phi)$ between two monoidal Doi Hom-Hopf datum is called \emph{monoidal} if the $\varphi$ and $\phi$ are monoidal Hom-bialgebra maps. We now consider the particular situation where $H=H'$ and $A=A'$. The following result is a generalization of \cite{CR95}.

{\bf Theorem 4.3.} Let $\xi=(id_H, id_A, \phi): (H,A,C)\rightarrow (H,A,C')$ be a monoidal morphism of monoidal Doi Hom-Hopf datum. Then
\begin{eqnarray}
G(C')=C.
\end{eqnarray}
Let $(M,\mu)\in{}_A\widetilde{\mathscr{H}}(\mathscr{M}_k)(H)^{C} $ be flat as a $k$-module, and take $(N,\nu)\in{}_A\widetilde{\mathscr{H}}(\mathscr{M}_k)(H)^{C'}$. If $(C,\g)$ is a monoidal Hom-Hopf algebra, then
\begin{eqnarray}
M\o G(N)\cong G(F(M)\o N) ~~in~~{}_A\widetilde{\mathscr{H}}(\mathscr{M}_k)(H)^{C}.
\end{eqnarray}
If $(C,\g)$ has a twisted antipode $\overline{S}$, then
\begin{eqnarray}
G(N) \o M \cong G( N\o F(M)) ~~in ~~{}_A\widetilde{\mathscr{H}}(\mathscr{M}_k)(H)^{C}.
\end{eqnarray}

{\bf Proof.} We know that $\varepsilon_{C'}\o id_C: C'\square_CC\rightarrow C$ is an isomorphism; the inverse map is $(\phi\o id_C)\Delta_C: C\rightarrow C'\square_CC$. It is clear that $\varepsilon_{C'}\o id_C$ is $(A,\b)$-linear and $(C,\g)$-colinear. and this prove Eq.(4.7).

Now we define the map
\begin{eqnarray*}
\Gamma: M\o G(N)=M\o (N\square_{C'}C),
G(F(M)\o N)=(F(M)\o N)\square_{C'}C,
\end{eqnarray*}
which is given by
\begin{eqnarray*}
\Gamma(m\o (n_i\o c_i))=(m_{[0]}\o n_i)\o m_{[1]}c_i.
\end{eqnarray*}
Recall that $F(M)=M$ as an $(A,\b)$-module, with $(C', \g')$-coaction given by
\begin{eqnarray*}
\rho_{F(M)}(m)=m_{[0]}\o \phi(m_{[1]}).
\end{eqnarray*}
(1) $\Gamma$ is well-defined, we have to show that
\begin{eqnarray*}
\Gamma(m_i\o (n_i\o c_i))\in (F(M)\o N)\square_C'C.
\end{eqnarray*}
This may be seen as follows: for any $m\in M$ and $n_i\square_{C'}c\in N\square_{C'}C$, we have
\begin{eqnarray*}
&&(\rho_{F(M)\o N}\o id_C)((m_{[0]}\o n_i)\o m_{[1]}c_i)\\
&=&(m_{[0][0]}\o n_{i[0]})\o \phi(m_{[0][1]})n_{i[1]}\o m_{[1]}c_i\\
&=&(\mu(m_{[0]})\o \nu(n_{i}))\o \phi(m_{[0][1]})\phi(c_{i(1)})\o \g^{-1}(m_{[1]}c_{i(2)})\\
&=&(m_{[0]}\o n_{i})\o[ \phi(m_{[0][1]})\phi(c_{i(1)})\o m_{[1]}c_{i(2)}]\\
&=& (id_{F(M)\o N}\o \rho_{C'})((m_{[0]}\o n_i)\o m_{[1]}c_i).
\end{eqnarray*}
(2) $\Gamma$ is $(A,\b)$-linear. Indeed, for any $a\in A, m\in M$ and $n_i\square_{C'}c\in N\square_{C'}C$, we have
\begin{eqnarray*}
&&\Gamma(a\c (m\o (n_i\o c_i)))\\
&=& \Gamma(a_{(1)}\c m\o (a_{(2)[0]}\c n_i\o a_{(2)[1]}\c c_i))\\
&=& (a_{(1)[0]} \c m_{[0]} \o a_{(2)[0]}\c n_i)\o (a_{(1)[1]} \c m_{[1]} )(a_{(2)[1]} \c c_i)\\
&=& (a_{[0](1)} \c m_{[0]} \o a_{[0](2)}\c n_i)\o a_{(1)} \c (m_{[1]} c_i)\\
&=&a_{[0]}\c (m_{[0]}\o n_i)\o a_{(1)} \c (m_{[1]} c_i)\\
&=&a \c \Gamma(m\o (n_i\o c_i)).
\end{eqnarray*}
(3) $\Gamma$ is $(C,\g)$-colinear. Indeed, for any $m\in M$ and $n_i\square_{C'}c\in N\square_{C'}C$, we have
\begin{eqnarray*}
&&\rho\circ \Gamma(m\o (n_i\o c_i))\\
&=&\rho((m_{[0]}\o n_i)\o m_{[1]}c_i)\\
&=& (\mu^{-1}(m_{[0]})\o \nu^{-1}(n_i)) \o  m_{[1](1)}c_{i(1)}\o \g(m_{[1](2)}c_{i(2)})\\
&=& (m_{[0]}\o \nu^{-1}(n_i)) \o  m_{[0][1]}c_{i(1)}\o m_{[1]}\g(c_{i(2)})\\
&=& (\Gamma\o id_C)(m_{[0]}\o (\nu^{-1}(n_i)\o c_{i(1)}))\o m_{[1]}\g(c_{i(2)})\\
&=& (\Gamma\o id_C)\circ \rho(m\o (n_i\o c_i)).
\end{eqnarray*}
Assume $(C, \g)$ has an antipode, define
\begin{eqnarray*}
&&\Psi: (F(M)\o N)\square_{C'}C\rightarrow M\o (N\square_{C'}C),\\
&&\Psi((m_i\o n_i)\o c_i)=\mu^{2}(m_{i[0]})\o (n_i\o S(m_{i[1]})\g^{-2}(c_i)).
\end{eqnarray*}
We have to show that $\Psi$ is well-defined. $(M,\mu)$ is flat, so $M\o (N\square_{C'}C)$ is the equalizer of the maps
\begin{eqnarray*}
id_M\o id_N\o \rho_C:~~~  M\o N\o C\rightarrow M\o N\o C'\o C
\end{eqnarray*}
and \begin{eqnarray*}
   id_M\o \rho_N\o id_C:~~  M\o N\o C\rightarrow M\o N\o C'\o C.
    \end{eqnarray*}
    Now take $(m_i\o n_i)\o c_i\in (F(M)\o N)\square_{C'}C$, then
    \begin{eqnarray}
   (m_{i[0]}\o n_{i[0]})\o \phi(m_{i[1]})n_{i[1]}\o c_i=(\mu^{-1}(m_i)\o \nu^{-1}(n_i))\o \phi(c_{i(1)})\o \g(c_{i(2)}).
    \end{eqnarray}
Therefore, we get
    \begin{eqnarray*}
   &&id_M\o id_N\o \rho_C(\mu^{2}(m_{i[0]})\o (n_i\o S(m_{i[1]})\g^{-2}(c_i)))\\
   &=&\mu^{2}(m_{i[0]})\o (n_i\o \phi(S(m_{i[1](2)})\g^{-2}(c_{i(1)}))\o S(m_{i[1](1)})\g^{-2}(c_{i(2)}))\\
    &=&m_{i[0]}\o \nu^{-1}(n_i)\o \phi(S(\g(m_{i[1](2)}))\g^{-1}(c_{i(1)}))\o S(\g^{2}(m_{i[1](1)}))c_{i(2)}
    \end{eqnarray*}
and
\begin{eqnarray*}
&&id_M\o \rho_N\o id_C(\mu^{2}(m_{i[0]})\o (n_i\o S(m_{i[1]})\g^{-2}(c_i)))\\
&=&\mu^{2}(m_{i[0]})\o (n_{i[0]}\o n_{i[1]}\o S(m_{i[1]})\g^{-2}(c_i))\\
&=& m_{i[0]}\o n_{i[0]}\o \g(n_{i[1]})\o S(\g(m_{i[1]}))\g^{-1}(c_i).
\end{eqnarray*}
Applying $(id_M\o \phi\o id_C)\circ(id_M\o (\Delta_C\circ S_C))\circ\rho_M$ to the first factor of Eq.(4.10), we obtain
\begin{eqnarray*}
&&m_{i[0][0]}\o \phi(S(m_{i[0][1](2)}))\o S(m_{i[0][1](1)}) \o n_{i[0]}\o \phi(m_{i[1]})n_{i[1]}\o c_i\\
&=& \mu^{-1}(m_{i[0]})\o \phi(S(\g^{-1}(m_{i[1](2)})))\o S(\g^{-1}(m_{i[1](1)})) \o \nu^{-1}(n_{i})\o \phi(c_{i(1)})\o \g(c_{i(2)}).
\end{eqnarray*}
Applying $id_M \o \g^{2}\o id_C\o id_N \o \g^{-1} \o  \g^{-1} $ to the above identity, we have
\begin{eqnarray*}
&&m_{i[0][0]}\o \phi(S(\g^{2}(m_{i[0][1](2)})))\o S(m_{i[0][1](1)}) \o n_{i[0]}\o \g^{-1}(\phi(m_{i[1]})n_{i[1]})\o \g^{-1}(c_i)\\
&=& \mu^{-1}(m_{i[0]})\o \phi(S(\g(m_{i[1](2)})))\o S(\g^{-1}(m_{i[1](1)})) \o \nu^{-1}(n_{i})\o \phi(\g^{-1}(c_{i(1)}))\o c_{i(2)}.
\end{eqnarray*}
Multiplying the second and the fifth  factor, and also the third and sixth factor, we have
\begin{eqnarray*}
&&\mu(m_{i[0]})\o n_{i[0]}\o \g(n_{i[1]})\o S(\g(m_{i[1]}))\g^{-1}(c_i)\\
&=&\mu(m_{i[0]})\o \nu^{-1}(n_i)\o \phi(S(\g(m_{i[1](2)}))\g^{-1}(c_{i(1)}))\o S(\g^{2}(m_{i[1](1)}))c_{i(2)},
\end{eqnarray*}
and applying $\mu^{-1} \o id_N \o id_C \o id_C $ to the above identity, we obtain
\begin{eqnarray*}
&&m_{i[0]}\o n_{i[0]}\o \g(n_{i[1]})\o S(\g(m_{i[1]}))\g^{-1}(c_i)\\
&=&m_{i[0]}\o \nu^{-1}(n_i)\o \phi(S(\g(m_{i[1](2)}))\g^{-1}(c_{i(1)}))\o S(\g^{2}(m_{i[1](1)}))c_{i(2)}
\end{eqnarray*}
or
\begin{eqnarray*}
id_M\o \rho_N\o id_C\circ(\Psi((m_i\o n_i)\o c_i))
= id_M\o id_N\o \rho_C\circ(\Psi((m_i\o n_i)\o c_i)).
\end{eqnarray*}
Let us point out that $\Gamma$ and $\Psi$ are each other's inverses. In fact,
\begin{eqnarray*}
&&\Gamma\circ\Psi((m_i\o n_i)\o c_i)\\
&=&\Gamma(\mu^{2}(m_{i[0]})\o (n_i\o S(m_{i[1]}\g^{-2}(c_i))))\\
&=& (\mu^{2}(m_{i[0][0]})\o n_i)\o  \g^{2}(m_{i[0][1]})S(m_{i[1]})\g^{-2}(c_i))\\
&=&(\mu^{2}(m_{i[0][0]})\o n_i)\o  [\g(m_{i[0][1]})S(m_{i[1]})]\g^{-1}(c_i))\\
&=&(\mu(m_{i[0]})\o n_i)\o  [\g(m_{i[1](1)})S(\g(m_{i[1](2)}))]\g^{-1}(c_i))\\
&=&(m_i\o n_i)\o c_i,
\end{eqnarray*}
and
\begin{eqnarray*}
&&\Psi\circ\Gamma(m\o (n_i\o c_i))\\
&=&\Psi((m_{[0]}\o n_i)\o m_{[1]}c_i)\\
&=&\mu^{2}(m_{[0][0]})\o (n_i\o [S(\g^{-1}(m_{[0][1]}))\g^{-2}(m_{[1]})]\g^{-1}(c_i))\\
&=&\mu(m_{[0]})\o (n_i\o [S(\g^{-1}(m_{[1](1)}))\g^{-1}(m_{[1](2)})]\g^{-1}(c_i))\\
&=&m\o (n_i\o c_i).
\end{eqnarray*}
The proof of Eq.4.9 is similar and  left to the reader. $\hfill \Box$

{\bf Corollary 4.4.} Let $(H,A,C)$ be a monoidal Doi Hom-Hopf Datum,  $\Lambda$: ${}_A\widetilde{\mathscr{H}}$$(\mathscr{M}_k)$$(H)^{C}$ $\rightarrow $ ${}_A\widetilde{\mathscr{H}}$$(\mathscr{M}_k)(H)$ the functor forgetting the $(C,\g)$-coaction. For any flat Doi Hom-Hopf module $(M, \mu)$, we have an isomorphism \begin{eqnarray*}
           M\o C\cong \Lambda(M)\o C
            \end{eqnarray*}
in ${}_A\widetilde{\mathscr{H}}(\mathscr{M}_k)(H)^{C}$. If $k$ is a field, then ${}_A\widetilde{\mathscr{H}}(\mathscr{M}_k)(H)^{C}$ has enough injective objects, and any injective object in  ${}_A\widetilde{\mathscr{H}}(\mathscr{M}_k)(H)^{C}$ is a direct summand of an object of the form $I\o C$, where $I$ is  an injective $(A,\b)$-module.

We have already  proved that the category of Hom-Yetter-Drinfeld modules may be viewed as the category of Doi Hom-Hopf modules corresponding to a monoidal Doi Hom-Hopf Datum. Then we have the following corollary.

{\bf Corollary 4.5.} Let $(H,\a)$ be a monoidal Hom-Hopf algebra with the bijective antipode. Then the category of Hom-Yetter-Drinfeld modules over $(H,\a)$ is a Grothendieck category with enough injective objects.

We continue with the dual version of Theorem 4.3.

{\bf Theorem 4.6.} Let $\chi=(id_H, \psi, id_C): (H,A,C)\rightarrow (H,A',C)$ be a monoidal morphism of monoidal Doi Hom-Hopf data. Then
\begin{eqnarray}
F(A)=A'.
\end{eqnarray}
Let $(M,\mu)\in{}_A\widetilde{\mathscr{H}}(\mathscr{M}_k)(H)^{C} $ be flat as a $k$-module, and take $(N,\nu)\in{}_{A'}\widetilde{\mathscr{H}}(\mathscr{M}_k)(H)^{C}$. If $(A', \b')$ is a monoidal Hom-Hopf algebra, then
\begin{eqnarray}
F(M)\o N\cong F(M\o G(N)) ~~in~~{}_A\widetilde{\mathscr{H}}(\mathscr{M}_k)(H)^{C}.
\end{eqnarray}
If $(A', \b')$ has a twisted antipode $\overline{S}$, then
\begin{eqnarray}
N \o F(M) \cong F( G(N)\o M) ~~in ~~{}_A\widetilde{\mathscr{H}}(\mathscr{M}_k)(H)^{C}.
\end{eqnarray}

{\bf Proof.} We only prove Eq.(4.12) and similar for  Eq.(4.11) and Eq.(4.13). Assume that $(A', \b')$ is a monoidal Hom-Hopf algebra and define
\begin{eqnarray*}
\Gamma: F(M\o G(N))=A'\o _A M\o G(N)\rightarrow F(M)\o N =(A'\o _A M)\o N
\end{eqnarray*}
 by
\begin{eqnarray*}
\Gamma(a'\o (m\o n))=(a'_{(1)}\o m) \o a'_{(2)}\c n,
\end{eqnarray*}
for all $a'\in A', m\in M$ and $n\in N$. $\Gamma$ is well-defined since
\begin{eqnarray*}
\Gamma(a'\psi(a)\o (m\o n))&=& (a'_{(1)}\psi(a_{(1)})\o m) \o a'_{(2)}\psi(a_{(2)})\c n\\
&=& (a'_{(1)}\o a_{(1)}\c m) \o a'_{(2)}\psi(a_{(2)})\c n\\
&=& \Gamma(a'\o (a_{(1)}\c m\o\psi(a_{(2)})\c n ))\\
&=& \Gamma(a'\o a\c ( m\o n )).
\end{eqnarray*}
It is easy to check that $\Gamma$ is $(A', \b')$-linear.
 Now we shall verify that $\Gamma$ is $(C,\g)$-colinear based on Eq.(3.1). For any $a'\in A', m\in M$ and $n\in N$, we have
\begin{eqnarray*}
\rho(\Gamma(a'\o (m\o n)))&=& \rho((a'_{(1)}\o m) \o a'_{(2)}\c n)\\
&=&(a'_{(1)[0]}\o m_{[0]}) \o (a'_{(2)[0]}\c n_{[0]})\o (a'_{(1)[1]}\o m_{[1]})(a'_{(2)[1]}\c n_{[1]})\\
&\stackrel{(3.1)}{=}& (a'_{[0](1)}\o m_{[0]}) \o (a'_{[0](2)}\c n_{[0]})\o a'_{[1]}( m_{[1]}n_{[1]})\\
&=&(\Gamma\o id_c)(a'_{[0]}\o ( m_{[0]}\o  n_{[0]}))\o a'_{[1]}( m_{[1]}n_{[1]})\\
&=&(\Gamma\o id_c)\rho(a'\o (m\o n)).
\end{eqnarray*}
The inverse of $\Gamma$ is given by
\begin{eqnarray*}
\Psi((a'\o m)\o n)=\b'^{2}(a'_{(1)})\o (m\o S(a'_{(2)})\nu^{-2}(n))
\end{eqnarray*}
for all $a'\in A', m\in M$ and $n\in N$. One can check that $\Psi$ is well-defined  similar to $\Gamma$.
Finally, we have
 \begin{eqnarray*}
        \Psi(\Gamma(a'\o (m\o n)))&=& \Psi((a'_{(1)}\o m) \o a'_{(2)}\c n)\\
        &=&\b'^{2}(a'_{(1)(1)})\o (m\o S(a'_{(1)(2)})\nu^{-2}(a'_{(2)}\c n))\\
        &=&\b'(a'_{(1)})\o (m\o [S(\b'^{-1}(a'_{(2)(1)}))\b'^{-1}(a'_{(2)(2)}]\c \nu^{-1}(n))\\
        &=&a'\o a'\o (m\o n)
         \end{eqnarray*}
         and
         \begin{eqnarray*}
        \Gamma(\Psi((a'\o m)\o n))&=&\Gamma(\b'^{2}(a'_{(1)})\o (m\o S(a'_{(2)})\nu^{-2}(n)))\\
        &=&(\b'^{2}(a'_{(1)(1)})\o m) \o a'_{(2)}\c \b'^{2}(a'_{(1)(2)})\c S(a'_{(2)})\nu^{-2}(n)\\
        &=&(\b'(a'_{(1)})\o m) \o a'_{(2)}\c [\b'(a'_{(2)(1)})\c S(\b'(a'_{(2)}))]\nu^{-1}(n)\\
        &=& (a'\o m)\o n),
         \end{eqnarray*}
as needed. The proof is completed. $\hfill \Box$

\section{Braidings on the category of Doi Hom-Hopf modules}
\def\theequation{5. \arabic{equation}}
\setcounter{equation} {0}

Consider now a map $\mathscr{Q}: \,\, C\o C\rightarrow A\o A$, with
 twisted convolution inverse $\mathscr{R}$ such that $(\b\o \b)\mathscr{Q}=\mathscr{Q}(\g\o \g)$ and $(\b\o \b)\mathscr{R}=\mathscr{R}(\g\o \g)$. This means that
\begin{eqnarray*}
 &&\mathscr{R}(\mathscr{Q}^1(c_{(2)}\o d_{(2)})_{[1]} \c \gamma^{-1}(c_{(1)})  \o \mathscr{Q}^2(c_{(2)}\o d_{(2)})_{[1]} \c \gamma^{-1}(d_{(1)}))(\beta(\mathscr{Q}^2(c_{(2)}\o d_{(2)})_{[0]}) \\
 &&~~~~\o \beta(\mathscr{Q}^1(c_{(2)}\o d_{(2)})_{[0]}))\\
&&~~~~~~~~ =\varepsilon_C(c)1_A \o \varepsilon_C(d)1_A,
\end{eqnarray*}
\begin{eqnarray*}
 &&\mathscr{Q}(\mathscr{R}^2(c_{(2)}\o d_{(2)})_{[1]} \c \gamma^{-1}(c_{(1)})  \o \mathscr{R}^1(c_{(2)}\o d_{(2)})_{[1]} \c \gamma^{-1}(d_{(1)}))(\beta(\mathscr{R}^2(c_{(2)}\o d_{(2)})_{[0]})\\
 &&~~~~ \o \beta(\mathscr{R}^1(c_{(2)}\o d_{(2)})_{[0]}))\\
&&~~~~~~~~ =\varepsilon_C(c)1_A \o \varepsilon_C(d)1_A,
\end{eqnarray*}
 for all $c, d\in C$. Sometimes, we write
 $\mathscr{Q}(c\o d):=\mathscr{Q}^1(c\o d)\o \mathscr{Q}^2(c\o d)$
 for all $c, d\in C$.

 Let $(M, \mu), (N,\nu) \in {}_A\widetilde{\mathscr{H}}(\mathscr{M}_k)(H)^{C}$.
 By Proposition 3.3  we know $(M\o N, \mu\o \nu)\in {}_A\widetilde{\mathscr{H}}(\mathscr{M}_k)(H)^{C}$.
Define a map
\begin{eqnarray}
&&c_{M, N}: M \o N \rightarrow N \o M, \nonumber\\
&&c_{M, N}(m\o n)=\mathscr {Q}(n_{[1]}\o
m_{[1]})(n_{[0]}\o m_{[0]}).
\end{eqnarray}

Next we will prove  that $c_{M,N}$ is an isomorphism with an inverse map
\begin{eqnarray*}
&&c^{-1}_{M, N}: N \o M \rightarrow M \o N, \nonumber\\
&&c^{-1}_{M, N}(n\o m)=\mathscr{R}(n_{[1]}\o
m_{[1]})(m_{[0]}\o n_{[0]}).
\end{eqnarray*}
For any $m\in M$ and $n\in N$, we have
\begin{eqnarray*}
&&c^{-1}_{M, N}\circ c_{M,N}(m\o n)\\
&=& c^{-1}_{M, N} (\mathscr{Q}(n_{[1]}\o
m_{[1]})(n_{[0]}\o m_{[0]}))\\
&=& \mathscr{R}( (\mathscr{Q}^1(n_{[1]}\o m_{[1]}) \c n_{[0]})_{[1]} \o (\mathscr{Q}^2(n_{[1]}\o m_{[1]}) \c m_{[0]})_{[1]}) \\
&&~~~~
( (\mathscr{Q}^2(n_{[1]}\o m_{[1]}) \c m_{[0]})_{[0]} \o (\mathscr{Q}^1(n_{[1]}\o m_{[1]}) \c n_{[0]})_{[0]}) \\
&=& \mathscr{R}(\mathscr{Q}^1(\gamma(n_{[1](2)})\o \gamma(m_{[1](2)}))_{[1]} \c n_{[1](1)} \o \mathscr{Q}^2(\gamma(n_{[1](2)})\o \gamma( m_{[1](2)}))_{[1]} \c m_{[1](1)}) \\
&&~~~~
(\mathscr{Q}^2(\gamma(n_{[1](2)})\o \gamma(m_{[1](2)}))_{[0]} \c \mu^{-1}(m_{[0]}) \o \mathscr{Q}^1(\gamma(n_{[1](2)})\o \gamma( m_{[1](2)}))_{[0]} \c \nu^{-1}(n_{[0]}))\\
&=& (\mathscr{R}(\mathscr{Q}^1(n_{[1](2)} \o m_{[1](2)})_{[1]} \c\gamma^{-1}(n_{[1](1)})  \o \mathscr{Q}^2(n_{[1](2)} \o m_{[1](2)})_{[1]} \c\gamma^{-1}(m_{[1](1)})) \\
&&~~~~(\beta(\mathscr{Q}^2(n_{[1](2)} \o m_{[1](2)})_{[0]}) \o \beta(\mathscr{Q}^1(n_{[1](2)} \o m_{[1](2)})_{[0]})))(m_{[0]} \o n_{[0]})
\end{eqnarray*}
\begin{eqnarray*}
&=& (\varepsilon_C(m_{[1]})1_A \o \varepsilon_C(n_{[1]})1_A)(m_{[0]} \o n_{[0]}) \\
&=& m\o n.
\end{eqnarray*}
So $c^{-1}_{M, N}\circ c_{M,N}=id_{M\o N}$.
Similarly, we can prove $c_{M,N}\circ c^{-1}_{M, N}=id_{N\o M}$.

Our aim is now to give necessary and sufficient conditions on
$\mathscr{Q}$ such that the $c_{M, N}$ defines a braiding on the
 monoidal category of Doi Hom-Hopf modules. 
 Recall from \cite{LS14} that for any $(M, \mu), (N,\nu) \in {}_A\widetilde{\mathscr{H}}(\mathscr{M}_k)(H)^{C}$, 
 the associativity and unit constraints are given by
\begin{eqnarray*}
&&a_{M,N,P}:(M\o N)\o P\rightarrow M\o (N\o P),~(m\o n)\o p\mapsto \mu(m)\o (n\o \pi^{-1}(p)),\\
&& l_M: k\o M\rightarrow M,~~k\o m\mapsto k\mu(m), r_M: M\o k\rightarrow M, ~~m\o k\mapsto k\mu(m).
\end{eqnarray*}
 Next, we will find conditions under which $c_{M, N}$ is both an
 $(A, \b)$-module map and a $(C, \g)$-comodule map,
 and  satisfies the following conditions (for $P \in {}_A\widetilde{\mathscr{H}}(\mathscr{M}_k)(H)^{C}$)
\begin{equation}
 a_{N,P,M}\circ c_{M,    N\o P}\circ a_{M,N,P}=(id_{N} \o c_{M,P} )\circ a_{N,M,P}\circ(c_{M,N}\o id_{P}).
\end{equation}
\begin{equation}
a^{-1}_{N,P,M}\circ c_{M\o N, P}\circ a^{-1}_{M,N,P}=(c_{M, P}\o id_{N})\circ a^{-1}_{M,P,N} \circ(id_{M}\o
c_{N,P}),
\end{equation}

Recall from \cite{GC14} that $A\o C$ can be made into a Doi Hom-Hopf module as follows:
the $(A,\b)$-action and $(C,\g)$-coaction on $A\o C$ are given by the formulas
$$
\left\{
  \begin{array}{ll}
a\c (b \otimes c )   =\b^{-1}( a) b \otimes \g(c);\\
  \rho_{A\o C}  ( b \otimes c  )
  =(b_{[0]} \otimes c _{(1)})\otimes b_{[1]}c_{(2)},
\end{array}
\right.
$$
for any $a, b\in A$ and $c\in C$.

In order to approach to our main result, we need some lemmas:

{\bf Lemma 5.1.} Let $(M, \mu), (N,\nu) \in {}_A\widetilde{\mathscr{H}}(\mathscr{M}_k)(H)^{C}$. Then
 $c_{M, N}$ is $(A, \b)$-linear if and only if the following condition
 is satisfied:
\begin{eqnarray}
&&\mathscr{Q}(a_{(2)[1]}
 \c c\o
a_{(1)[1]}\c d)
 (a_{(2)[0]}\o a_{(1)[0]})
=\D (a)\mathscr{Q}(c\o d)
\end{eqnarray}
for all $a\in A$ and $c, d\in C$.

{\bf Proof. } If  $c_{M, N}$ is $(A, \b)$-linear then $a\triangleright c_{M, N}(m\o n)=
 c_{M, N}(a\triangleright (m\o n))$. We compute two sides of the equation as
 follows:
  $$
 a\triangleright c_{M, N}(m\o n)
=(a_{(1)}\o a_{(2)})\mathscr {Q}(n_{[1]}\o
m_{[1]})(n_{[0]}\o m_{[0]})
$$
and
\begin{eqnarray*}
&&c_{M, N}(a\triangleright (m\o n))\\
&&=\mathscr{Q}(a_{(2)[1]}
 \c n_{[1]}\o
a_{(1)[1]}\c m_{[1]})
 (a_{(2)[0]}\c n_{[0]}\o a_{(1)[0]}\c m_{[0]}).
\end{eqnarray*}

 Conversely, considering these equations and taking $M=N=A\o C$ and  $m=1\o c$
 and $n=1\o d$ for all $c, d\in C$. Then we can get Eq.(5.4). $\hfill \Box$

{\bf Definition 5.2.} A \emph{quasitriangular monoidal Hom-Hopf algebra} is a monoidal Hom-Hopf algebra $(H, \alpha)$ together with an invertible element
$R=R^{(1)}\o R^{(2)}\in H\o H$  such that
he following conditions hold:
\begin{eqnarray*}
&&(QT1)~\Delta(R^{(1)})\otimes R^{(2)}=R^{(1)}\otimes r^{(1)}\otimes R^{(2)}r^{(2)},\\
&&(QT2)~R^{1}\otimes\Delta(R^{2})=R^{1}r^{1}\otimes r^{2}\otimes R^{2},\\
&&(QT3)~\varepsilon(R^{(1)})R^{(2)}=1_{H},R^{(1)}\varepsilon(R^{(2)})=1_{H},\\
&&(QT4)~\Delta^{cop}(h)R=R\Delta(h),\\
&&(QT5)~(\a \o \a)(R) = R,
\end{eqnarray*}
where $\Delta^{cop}(h)=h_{(2)}\otimes h_{(1)}$ for all $h\in H$.  $(H, \alpha)$ is called almost cocommutative  if $\Delta^{cop}(h)R=R\Delta(h)$ holds.

{\bf Example 5.3.} Suppose that $C=k$ and write $R=\mathscr{Q}(1\o 1)$. Then Eq.(5.4) takes the form $R\Delta_A^{cop}(a)=\Delta_A(a)R$ and this means that $(A,\b)$ is almost cocommutative.

{\bf Lemma 5.4.} Let $(M, \mu),(N,\nu) \in {}_A\widetilde{\mathscr{H}}(\mathscr{M}_k)(H)^{C}$. Then
 $c_{M, N}$ is $(C, \g)$-colinear if and only if the following condition
 is satisfied:
\begin{eqnarray}
&&\mathscr{Q}(d_{(2)}\o c_{(2)})_{[0]}\o m_C(\mathscr{Q}(d_{(2)}\o c_{(2)})_{[1]}(d_{(1)}\o c_{(1)}))\nonumber\\
 &=&\mathscr{Q}(d_{(1)}\o c_{(1)})\o c_{(2)}d_{(2)},
\end{eqnarray}
for all  $c, d\in C$.

{\bf Proof.}
 If $c_{M, N}$ is $(C,\g)$-colinear, then we do
 the following calculations:
\begin{eqnarray*}
&&\r_{N\o M}  c_{M, N}(m\o n)\\
&=&\r_{N\o M} (\mathscr {Q}(n_{[1]}\o
m_{[1]})(n_{[0]}\o m_{[0]}))\\
&=& \mathscr {Q}(n_{[1]}\o
m_{[1]})_{[0]}(n_{[0][0]}\o m_{[0][0]})\o  m_C(\mathscr {Q}(n_{[1]}\o
m_{[1]})_{[1]}(n_{[0][1]}\o m_{[0][1]}))\\
&=&\mathscr {Q}(\gamma^{-1}(n_{[1](2)})\o
\gamma^{-1}(m_{[1](2)}))_{[0]}(\nu(n_{[0]})\o \mu(m_{[0]}))\o  m_C(\mathscr {Q}(\gamma^{-1}(n_{[1](2)})\\
&&\hspace{6cm}\o
\gamma^{-1}(m_{[1](2)}))_{[1]}(n_{[1](1)}\o m_{[1](2)})).
\end{eqnarray*}
On the other hand, we have
\begin{eqnarray*}
(c_{M, N}\o id_C) \r_{M\o N} (m\o n)&=&\mathscr{Q}(n_{[0][1]}\o
 m_{[0][1]})(n_{[0][0]}\o m_{[0][0]}) \o (m_{[1]}n_{[1]})\\
 &=&\mathscr{Q}(n_{[1](1)}\o
 m_{[1](1)})(\nu(n_{[0]})\o \mu(m_{[0]})) \o \gamma^{-1}(m_{[1](2)}n_{[1](2)}).
\end{eqnarray*}

 Conversely, let $M=N=A\o C$ and take $m=1\o c$ and $n=1\o d$
 for all $c, d\in C$. Then we can get Eq.(5.5).  $\hfill \Box$

Dual to quasitriangular monoidal Hom-Hopf algebras, a \emph{coquasitriangular monoidal Hom-Hopf algebra} is a monoidal Hom-Hopf algebra $(H, \alpha)$
together with a bilinear form $\sigma$ on $(H, \alpha)$ (i.e. $\sigma\in$ Hom($H\otimes H,k$))
such that the following axioms hold:
\begin{eqnarray*}
&&(BR1)~\sigma( hg,l)\rangle=\sigma( h,l_{(2)})\sigma( g,l_{(1)}),\\
&&(BR2)~\sigma( h,gl)=\sigma( h_{(1)},g)\sigma(h_{(2)},l),\\
&&(BR3)~\sigma( h_{(1)},g_{(1)}) g_{(2)}h_{(2)}=h_{(1)}g_{(1)}\sigma( h_{(2)},g_{(2)}),\\
&&(BR4)~\sigma( 1_H,h)=\sigma( h,1_H)=\varepsilon(h),\\
&&(BR5)~\sigma(\a(h),\a(g)) = \sigma(h,g),
\end{eqnarray*}
for all $h,g,l\in H$. $(H, \alpha)$ is called almost commutative  if $\sigma( h_{(1)},g_{(1)}) g_{(2)}h_{(2)}=h_{(1)}g_{(1)}$
$\sigma( h_{(2)},g_{(2)})$ holds.

{\bf Example 5.5.} Suppose $A=k$. Then Eq.(5.5) takes the form
\begin{eqnarray*}
\mathscr{Q}( h_{(1)},g_{(1)}) g_{(2)}h_{(2)}=h_{(1)}g_{(1)}\mathscr{Q}( h_{(2)},g_{(2)})
\end{eqnarray*}
and this means that $(A,\b)$ almost commutative.

{\bf Lemma 5.6.} Let $(M, \mu), (N,\nu), (P,\pi) \in {}_A\widetilde{\mathscr{H}}(\mathscr{M}_k)(H)^{C}$.
  Then Eq.(5.2) holds if and only if the following condition
 is satisfied, with $\mathscr{U}=\mathscr{Q}$:
\begin{eqnarray}
&&\mathscr {Q}^{1}(e\o
\gamma(d_{(2)}))\o (\mathscr{U}^1(\g^{-1}(c)\o \mathscr {Q}^{2}(e\o
\gamma(d_{(2)}))_{[1]}d_{(1)})\nonumber\\
&&\o \mathscr{U}^2(\g^{-2}(c)\o \mathscr {Q}^{2}(e\o
c_{(2)})_{[1]}\gamma^{-1}(c_{(1)}))\mathscr {Q}^{2}(e\o\gamma(d_{(2)}))_{[0]}\nonumber\\
&=&\mathscr {Q}^1(e\gamma^{-1}(c)\o \gamma(m_{[1]}))_{(1)}\o \mathscr {Q}^1(e\gamma^{-1}(c)\o \gamma(c))_{(2)}\o \mathscr {Q}^2(\g^{-1}(e)\gamma^{-2}(c)\o d)~~~~
 \end{eqnarray}
for all  $c, d, e\in C$.

{\bf Proof.} If Eq.(5.2) holds, then we compute
 as follows:
\begin{eqnarray*}
&&(id_{N} \o c_{M,P} )\circ a_{N,M,P}\circ(c_{M,N}\o id_{P})((m\o n)\o p)\\
&=&(id_{N} \o c_{M,P} )\circ a_{N,M,P}(\mathscr {Q}^{1}(n_{[1]}\o
m_{[1]})n_{[0]}\o \mathscr {Q}^{2}(n_{[1]}\o
m_{[1]})m_{[0]}\o p)\\
&=&(id_{N} \o c_{M,P} )(\b(\mathscr {Q}^{1}(n_{[1]}\o
m_{[1]}))\nu(n_{[0]})\o (\mathscr {Q}^{2}(n_{[1]}\o
m_{[1]})m_{[0]}\o \pi^{-1}(p)))\\
&=&\b(\mathscr {Q}^{1}(n_{[1]}\o
m_{[1]}))\nu(n_{[0]})\o \mathscr{U}(\g^{-1}(p_{[1]})\o \mathscr {Q}^{2}(n_{[1]}\o
m_{[1]})_{[1]}m_{[0][1]})\\
&&\hspace{7cm}(\pi^{-1}(p_{[0]})\o \mathscr {Q}^{2}(n_{[1]}\o
m_{[1]})_{[0]}m_{[0][0]})\\
&=&\b(\mathscr {Q}^{1}(n_{[1]}\o
\gamma(m_{[1](2)})))\nu(n_{[0]})\o \mathscr{U}(\g^{-1}(p_{[1]})\o \mathscr {Q}^{2}(n_{[1]}\o
\gamma(m_{[1](2)}))_{[1]}m_{[1](1)})\\
&&\hspace{6cm}(\pi^{-1}(p_{[0]})\o \mathscr {Q}^{2}(n_{[1]}\o
\gamma(m_{[1](2)}))_{[0]}\mu^{-1}(m_{[0]}))\\
&=&\b(\mathscr {Q}^{1}(n_{[1]}\o
\gamma(m_{[1](2)})))\nu(n_{[0]})\o (\mathscr{U}^1(\g^{-1}(p_{[1]})\o \mathscr {Q}^{2}(n_{[1]}\o
\gamma(m_{[1](2)}))_{[1]}m_{[1](1)})\\
&&\pi^{-1}(p_{[0]})\o \b^{-1}(\mathscr{U}^2(\g^{-1}(p_{[1]})\o \mathscr {Q}^{2}(n_{[1]}\o
\gamma(m_{[1](2)}))_{[1]}m_{[1](1)}))\\
&&\hspace{8cm}\mathscr {Q}^{2}(n_{[1]}\o\gamma(m_{[1](2)}))_{[0]}m_{[0]}).
\end{eqnarray*}
Also we have
\begin{eqnarray*}
&& a_{N,P,M}\circ c_{M,    N\o P}\circ a_{M,N,P}((m\o n)\o p)\\
&=& a_{N,P,M}\circ c_{M,    N\o P}(\mu(m)\o( n\o \pi^{-1}(p)))\\
&=&a_{N,P,M}((\Delta_A\o id_A)(\mathscr {Q}(n_{[1]}\gamma^{-1}(p_{[1]})\o \gamma(m_{[1]})))((n_{[0]}\o \pi^{-1}(p_{[0]})) \o \mu(m_{[0]})))\\
&=&\b(\mathscr {Q}^1(n_{[1]}\gamma^{-1}(p_{[1]})\o \gamma(m_{[1]}))_{(1)})\nu(n_{[0]})\o \mathscr {Q}^1(n_{[1]}\gamma^{-1}(p_{[1]})\o \gamma(m_{[1]}))_{(2)}\pi^{-1}(p_{[0]})\\
&&\hspace{7cm}\o \b^{-1}(\mathscr {Q}^2(n_{[1]}\gamma^{-1}(p_{[1]})\o \gamma(m_{[1]})))m_{[0]}.
\end{eqnarray*}

 Conversely, take $M=N=P=A\o C$ and $m=1\o d$, and $n=1\o e$, and $p=1\o c$
 for all $c, d, e\in C$. Then we obtain
 Eq.(5.6).
  $\hfill \Box$

The  proof of the following lemma is similar to that of the
 above   lemma.

{\bf Lemma 5.7. } Let $(M, \mu), (N,\nu), (P,\pi) \in {}_A\widetilde{\mathscr{H}}(\mathscr{M}_k)(H)^{C}$.
  Then Eq.(5.3) holds if and only if the following condition
 is satisfied, with $\mathscr{U}=\mathscr{Q}$:
\begin{eqnarray}
&&\mathscr{U}^1(\mathscr{Q}^1(c_{(2)}\o \g^{-1}(e))_{[1]}\g^{-1}(c_{(1)})\o \g^{-2}(d))\mathscr{Q}^1(\g(c_{(2)})\o e)_{[0]}\nonumber\\ &&\hspace{4cm}\o\mathscr{U}^2(\mathscr{Q}^1(\g(c_{(2)})\o e)_{[1]}c_{(1)}\o \g^{-1}(d))\o  \mathscr{Q}^2(c\o e)\nonumber \\
 && =\mathscr{Q}^1(c\o
\g^{-2}( d)\g^{-1}(e))\o \mathscr{Q}^2(\g(c)\o
\g^{-1}( d)e)_{(1)}\o \mathscr{Q}^2(\g(c)\o
\g^{-1}( d)e)_{(2)}~~
\end{eqnarray}
for all  $c, d, e\in C$.

{\bf Proof.} If Eq.(5.3) holds, it follows that
\begin{eqnarray*}
&&(c_{M, P}\o id_{N})\circ a^{-1}_{M,P,N} \circ(id_{M}\o
c_{N,P})(m\o (n\o p))\\
&=&(c_{M,  P}\o id_{N})\circ a^{-1}_{M,P,N} (m\o \mathscr{Q}(p_{[1]}\o n_{[1]})(p_{[0]}\o n_{[0]})) \\
&=&(c_{M,  P}\o id_{N}) ((\mu^{-1}(m)\o \mathscr{Q}^1(p_{[1]}\o n_{[1]})p_{[0]})\o \b(\mathscr{Q}^2(p_{[1]}\o n_{[1]}))\nu(n_{[0]})) \\
&=& \mathscr{U}(\mathscr{Q}^1(p_{[1]}\o n_{[1]})_{[1]}p_{[0][1]}\o \g^{-1}(m_{[1]}))(\mathscr{Q}^1(p_{[1]}\o n_{[1]})_{[0]}p_{[0][0]}\o \mu^{-1}(m_{[0]}))\\
&&\hspace{8cm}\o  \b(\mathscr{Q}^2(p_{[1]}\o n_{[1]}))\nu(n_{[0]}) \\
&=& \{\b^{-1}(\mathscr{U}^1(\mathscr{Q}^1(\g(p_{[1](2)})\o n_{[1]})_{[1]}p_{[1](1)}\o \g^{-1}(m_{[1]})))\mathscr{Q}^1(\g(p_{[1](2)})\o n_{[1]})_{[0]}\}p_{[0]}\\ &&\o\mathscr{U}^2(\mathscr{Q}^1(\g(p_{[1](2)})\o n_{[1]})_{[1]}p_{[1](1)}\o \g^{-1}(m_{[1]}))\mu^{-1}(m_{[0]})\o  \b(\mathscr{Q}^2(p_{[1]}\o n_{[1]}))\nu(n_{[0]})
\end{eqnarray*}
and
\begin{eqnarray*}
&&a^{-1}_{P,M,N}\circ c_{M\o N, P}\circ a^{-1}_{M,N,P}(m\o (n\o p))\\
&=&a^{-1}_{P,M,N}\circ c_{M\o N, P}((\mu^{-1}(m)\o n)\o \pi(p))\\
&=& a^{-1}_{P,M,N}\mathscr{Q}(\g(p_{[1]})\o
\g^{-1}( m_{[1]})n_{[1]})(\pi(p_{[0]})\o (\mu^{-1}(m_{[0]})\o n_{[0]}))\\
&=& \b^{-1}(\mathscr{Q}^1(\g(p_{[1]})\o
\g^{-1}( m_{[1]})n_{[1]}))p_{[0]}\o \mathscr{Q}^2(\g(p_{[1]})\o
\g^{-1}( m_{[1]})n_{[1]})_{(1)}\mu^{-1}(m_{[0]})\\
&&\hspace{6cm}\o \b(\mathscr{Q}^2(\g(p_{[1]})\o
\g^{-1}( m_{[1]})n_{[1]}))_{(2)}\nu(n_{[0]}).
\end{eqnarray*}

Conversely, take $M=N=P=A\o C$ and $m=1\o d$, and $n=1\o e$, and $p=1\o c$
 for all $c, d, e\in C$. Then we obtain
Eq.(5.7).
 $\hfill \Box$

 Therefore, we can summarize our results as follows.

{\bf Theorem 5.8.} Let  $(H, A, C)$ be a monoidal Doi Hom-Hopf datum,
 and $\mathscr{Q}: C\o C\lr A\o A$ a twisted convolution
 invertible map. For $(M, \mu), (N,\nu)\in {}_A\widetilde{\mathscr{H}}(\mathscr{M}_k)(H)^{C}$, then the family of maps
 $$
c_{M, N}: M \o N \rightarrow N \o M,\quad c_{M, N}(m\o
 n)=\mathscr {Q}(n_{[1]}\o m_{[1]})(n_{[0]}\o m_{[0]})
 $$
defines a braiding on the category of Doi Hom-Hopf modules ${}_A\widetilde{\mathscr{H}}(\mathscr{M}_k)(H)^{C}$  if and only
  Eqs.(5.4), (5.5), (5.6), (5.7) are satisfied.

{\bf Example 5.9.} (1) Take $A=k$ and write
\begin{eqnarray*}
R=\mathscr {Q}(1_C\o 1_C)=\sum R^{(1)}\o R^{(2)}=\sum r^{(1)}\o r^{(2)}.
\end{eqnarray*}
Eqs.(5.6) and (5.7) take the form
\begin{eqnarray*}
&&~\Delta(R^{(1)})\otimes R^{(2)}=R^{(1)}\otimes r^{(1)}\otimes r^{(2)}R^{(2)},\\
&&~R^{(1)}\otimes\Delta(R^{(2)})=r^{(1)}R^{(1)}\otimes r^{(2)}\otimes R^{(2)}.
\end{eqnarray*}
and the braiding is
\begin{eqnarray*}
c_{M, N}: M \o N \rightarrow N \o M,
c_{M, N}(m\o n)=R^{(2)}\c \nu^{-1}(n)\o R^{(1)}\c\mu^{-1}(m).
\end{eqnarray*}
Assume that  $R$ is $\alpha$-invariant
(i.e. $\alpha(R^{(1)})\o\alpha(R^{(2)})=R^{(1)}\o R^{(2)}$). We conclude that the conditions of Theorem 5.8 are satisfied if and only if $(C, R^{-1})$  is a quasitriangular monoidal Hom-bialgebra.

(2) If $C=k$, then Eqs.(5.6) and (5.7) take the form
\begin{eqnarray*}
&&~\sigma( hg,l)\rangle=\sigma( h,l_{(1)})\sigma( g,l_{(2)}),\\
&&~\sigma( h,gl)=\sigma( h_{(1)},l)\sigma(h_{(2)},g).
\end{eqnarray*}
and the braiding is
\begin{eqnarray*}
&&c_{M, N}: M \o N \rightarrow N \o M, \nonumber\\
&&c_{M, N}(m\o n)= \sigma(n_{[1]}, m_{[1]})\nu(n_{[0]})\o \mu(m_{[0]}).
\end{eqnarray*}
Assume that $\sigma$ is $\alpha$-invariant
 (i.e. $\sigma( \alpha(h),\alpha(g))=\sigma( h, g)$ for all $h,g\in H$)
and we conclude that the conditions of Theorem 5.8 are satisfied if and only if $(A, \sigma)$  is a coquasitriangular monoidal Hom-bialgebra.

(3) Let $(H,\a)$ be a monoidal Hom-Hopf algebra with bijective antipode. We have seen that the category of Doi Hom-Hopf modules $_H\widetilde{\mathscr{H}}(\mathscr{M}_k)(H^{op} \o H)^{H}$ and the category of Hom-Yetter-Drinfeld modules $_H\mathscr{HYD}^H$  are isomorphic. Recall from \cite{LS14} that $_H\mathscr{HYD}^H$ is a braided category. The braided is induced by
\begin{eqnarray*}
c_{M,N}: M\o N\rightarrow N\o M, ~~~m\o n\mapsto \nu(n_{[0]})\o n_{[1]}\mu^{-1}(m).
\end{eqnarray*}
The corresponding map $\mathscr {Q}$ is
\begin{eqnarray*}
\mathscr {Q}: H\o H\rightarrow H\o H, ~~~~h\o k\mapsto \eta(\varepsilon(k))\o h.
\end{eqnarray*}
It is straightforward to check that $\mathscr {Q}$ satisfies the conditions of Theorem 5.8.  $\hfill \Box$

\section{The smash product of Hom-bialgebras and the Drinfel'd double}
\def\theequation{6. \arabic{equation}}
\setcounter{equation} {0}

Let $(A,\b)$ be a right $(H,\a)$-Hom comodule algebra, and $(B, \zeta )$ a left $(H,\a)$-Hom module coalgebra.
Consider the following version of smash product $A\# B$ with the multiplication given by
\begin{eqnarray*}
(a\# b)(c\# d)=a\b(c_{[0]})\# (\zeta^{-1}(b)\leftharpoonup c_{[1]})d.
\end{eqnarray*}
Then $A\# B$ is a Hom associative algebra with unit $1_A\# 1_B$.

If $(C, \g)$ is a faithfully projective left  $(H,\a)$-Hom module coalgebra, Then $(C^{\ast}, \g^{\ast})$ is a right $(H, \a)$-Hom module algebra. The right $(H, \a)$-action is given by
\begin{eqnarray*}
(c^{\ast}\leftharpoonup h, c)=(c^{\ast}, h\c c).
\end{eqnarray*}
Given $(M,\mu) \in {}_A\widetilde{\mathscr{H}}(\mathscr{M}_k)(H)^{C}$, we define an $A\# C^{\ast}$-action on $(M,\mu)$ as follows
\begin{eqnarray*}
(a\# c^{\ast})\c m=<c^{\ast}, m_{[1]}>a\c m_{[0]}.
\end{eqnarray*}
Assume that $(A,\b)$ and $(B, \zeta )$ are both monoidal Hom-bialgebras, consider $\D_{A\#B}$ and $\varepsilon_{A\#B}$ defined by
\begin{eqnarray*}
&&\D_{A\#B}(a\#b)=(a_{(1)}\#b_{(1)})\o (a_{(2)}\#b_{(2)}),\\
&&\varepsilon_{A\#B}(a\#b)=\varepsilon_A(a)\varepsilon_B(b).
\end{eqnarray*}

{\bf Proposition 6.1.} With  notations as above. If
\begin{eqnarray}
&&\D_A(\b(a_{[0]}))\o \D_A(\zeta^{-1}(b)\leftharpoonup a_{[1]})\nonumber\\
&=& \b(a_{(1)[0]})\o \b(a_{(2)[0]})\o (\zeta^{-1}(b_{(1)})\leftharpoonup a_{[1](1)}) \o (\zeta^{-1}(b_{(2)})\leftharpoonup a_{[1](2)})
\end{eqnarray}
and
\begin{eqnarray}
\varepsilon_A(a_{[0]})\o \varepsilon_B(b\leftharpoonup a_{[1]})=\varepsilon_A(a)\varepsilon_B(b),
\end{eqnarray}
for all $a\in A$ and $b\in B$, then $A\#B$ is a  monoidal
Hom-bialgebra. If  $(A,\b)$ and $(B, \zeta )$ are both monodial Hom-Hopf algebras, then $A\#B$ is a monoidal Hom-Hopf algebras with the antipode given by
\begin{eqnarray*}
S_{A\#B}(a\#b)=S(\b(a))_{[0]}\# (S(\zeta^{-1}(b))\leftharpoonup S(a)_{[1]}).
\end{eqnarray*}

{\bf Proof.}  We leave it to the reader to show that $\D_{A\#B}$ is multiplicative if and only if Eq.(6.1) holds, and  $\varepsilon_{A\#B}$ is multiplicative if and only if Eq.(6.2) holds. We show that the antipode defined above is convolution invertible.  In fact, we have
\begin{eqnarray*}
&&(a_{(1)}\#b_{(1)})S_{A\#B}(a_{(2)}\#b_{(2)})\\
&=&(a_{(1)}\#b_{(1)})(S(\b(a_{(2)}))_{[0]}\o (S(\zeta^{-1}(b_{(2)}))\leftharpoonup S(a_{(2)})_{[1]}))\\
&=&a_{(1)}S(\b^{2}(a_{(2)}))_{[0][0]}\#(\zeta^{-1}(b_{(1)})\leftharpoonup S(\b(a_{(2)}))_{[0][1]})(S(\zeta^{-1}(b_{(2)}))\leftharpoonup S(a_{(2)})_{[1]}))\\
&=&a_{(1)}S(\b(a_{(2)}))_{[0]}\#(\zeta^{-1}(b_{(1)})\leftharpoonup S(\b(a_{(2)}))_{[1](1)})(S(\zeta^{-1}(b_{(2)}))\leftharpoonup S(\b(a_{(2)}))_{[1](2)}))\\
&=&a_{(1)}S(\b(a_{(2)}))_{[0]}\#(\zeta^{-1}(b_{(1)})S(\zeta^{-1}(b_{(2)})))\leftharpoonup S(\b(a_{(2)}))_{[1]}\\
&=& \varepsilon_A(a)\varepsilon_B(b),
\end{eqnarray*}
and
\begin{eqnarray*}
&&S_{A\#B}(a_{(1)}\#b_{(1)})(a_{(2)}\#b_{(2)})\\
&=&(S(\b(a_{(1)}))_{[0]}\o (S(\zeta^{-1}(b_{(1)}))\leftharpoonup S(a_{(1)})_{[1]}))(a_{(2)}\#b_{(2)})\\
&=&S(\b(a_{(1)}))_{[0]}\b(a_{(2)[0]})\# (S(\zeta^{-1}(b_{(1)}))\leftharpoonup S(a_{(1)})_{[1]}a_{(2)[1]})b_{(2)}\\
&=&\varepsilon_A(a)\varepsilon_B(b),
\end{eqnarray*}
as desired. $\hfill \Box$

{\bf Proposition 6.2.} Let $(H,A,C)$ be a monodial Doi Hom-Hopf module. Assume $(C,\g)$ is faithfully projective as a $k$-module. 
Then $(A,\b)$ and $(C^{\ast}, \g^{\ast})$ satisfy Eqs.(6.1), (6.2), ${}_A\widetilde{\mathscr{H}}(\mathscr{M}_k)(H)^{C}$ and $A\# C^{\ast}$-Hom modules as monoidal categories.

{\bf Proof.} We leave the proof  to the reader similar to \cite{CVZ98}.  $\hfill \Box$

 {\bf Example 6.3.} Assume that $(H,\a)$ is faithfully projective as a $k$-module. The monoidal Hom-algebra $A\#C^\ast$ is  nothing else than the Drinfel'd double $D(H)=H\# H^{\ast}$. 
 Then we can give the multiplication as  follows
 \begin{eqnarray*}
 (h\# f)(k\#g)=h\a^{2}(h_{(2)(1)})\# <\a^{\ast-2}(f), \a(h_{(2)(2)})\rightharpoonup \bullet \leftharpoonup S^{-1}(\a^{-1}(h))>g.
 \end{eqnarray*}

Now let   $(H, A, C)$ be a monoidal Doi Hom-Hopf datum,
 and $\mathscr{Q}: C\o C\lr A\o A$ a twisted convolution
 invertible map satisfying Eqs.(5.4), (5.5), (5.6), (5.7). $\mathscr{Q}$ induces the map
 \begin{eqnarray*}
\mathscr{\widetilde{Q}}: k\rightarrow (A\# C^{\ast})\o (A\# C^{\ast}).
 \end{eqnarray*}
 The braiding on ${}_A\widetilde{\mathscr{H}}(\mathscr{M}_k)(H)^{C}$ translates into a braiding on $A\# C^{\ast}$-Hom modules. This means that $A\# C^{\ast}$ is a quasitriangular monoidal Hom-Hopf algebra.  $\hfill \Box$
\section*{Acknowledgements}

 The work is supported by the NSF of Jiangsu Province (No. BK2012736),
   the Fund of Science and Technology Department of Guizhou Province (No. 2014GZ81365),
   the Anhui Province Excellent Young Talents Fund Project (No. 2013SQRL092ZD),
the Anhui Provincial Natural Science Foundation (No. 1408085QA06, 1408085QA08).


\begin{thebibliography}{00}


\bibitem{CG11}
S. Caenepeel and I. Goyvaerts,  Monoidal Hom-Hopf algebras.   {\it Comm. Algebra } {\bf 39} (2011) 2216-2240.

\bibitem{CM97}  S. Caenepeel,  G. Militaru and S. L. Zhu,    Crossed modules and Doi-Hopf modules.  {\it Israel J. Math.}  {\bf 100} (1997) 221-247.

\bibitem{CR95}
S. Caenepeel and S. Raianu, Induction functors for the Doi-Koppinen unified Hopf modules. in {\it Abelian Groups and Modules}, pp. 73-94, Kluwer Academic. Dordrecht, 1995.

\bibitem{CVZ98}
S. Caenepeel,  F. Van Oystaeyen and B. Zhou,  Making the category
 of Doi-Hopf modules into a braided monoidal category. {\it Algebras
 Representation Theory}.   {\bf 1} (1998) 75-96.



\bibitem{CZ13} Y. Y. Chen,  Z. W. Wang and L. Y. Zhang,  Integrals for monoidal Hom-Hopf algebras and their applications.
             {\em   J. Math. Phy.}  {\bf 54}(7) (2013) 073515.


\bibitem{CWZ13} Y. Y. Chen,  Z. W. Wang and L. Y. Zhang,   The FRT-type theorem for the Hom-Long equation. {\it Comm. Algebra } {\bf 41} (2013) 3931-3948.

\bibitem{CZ14} Y. Y. Chen and L. Y. Zhang. The category of Yetter-Drinfel'd Hom-modules and the
quantum Hom-Yang-Baxter equation. {\it J. Math. Phys.}  {\bf 55}(3) (2014) 031702.



\bibitem{D83}
Y. Doi, Unifying Hopf modules. {\it  J. Algebra } {\bf 153} (1992) 373-385.


\bibitem{FG10}
 Y. Fr\'{e}gier and A. Gohr,   On Hom type algebras. {\it J. Gen. Lie Theory Appl.} {\bf 4}(2010) 16.


\bibitem{FGS}
Y. Fr\'{e}gier, A. Gohr and S. D. Silvestrov,   Unital algebras of Hom-associative type
and surjective or injective twistings. {\it J. Gen. Lie Theory Appl.} {\bf 3} (2009) 285-295.



\bibitem{G10}
A. Gohr, On Hom-algebras with surjective twisting. {\it J. Algebra} {\bf 324} (2010) 1483-1491.

\bibitem{GC14}

S. J. Guo and X. L. Chen, A Maschke type theorem for relative Hom-Hopf modules. Accepted by {\it Czech. Math. J.} 2014.

\bibitem{HLS}
J. T. Hartwig, D. Larsson and S. D. Silvestrov,   Deformations of Lie algebras
using $\sigma$-derivations. {\it J. Algebra} {\bf 295} (2006) 314-361.


\bibitem{LS14} L. Liu and B. L. Sheng. Radford's biproducts and Yetter-Drinfeld modules for monoidal Hom-Hopf algebras.
{\em J. Math. Phys.}  {\bf 55}(3) (2014) 031701.

\bibitem{MP14} A. Makhlouf and F. Panaite. Yetter-Drinfeld modules for Hom-bialgebras.
{\em J. Math. Phys.}  {\bf 55}(1) (2014) 013501.

\bibitem{MS08}
 A. Makhlouf and  S. D. Silvestrov,   Hom-algebra stuctures. {\it J. Gen. Lie Theory Appl.}
{\bf 2} (2008) 51-64.

\bibitem{MS09}   A. Makhlouf, S. Silvestrov. Hom-Lie admissible Hom-coalgebras and Hom-Hopf algebras.
{\it J. Gen. Lie Theory in Mathematics, Physics and beyond}. Springer-Verlag,
Berlin, 2009, pp. 189-206.

\bibitem{MS10}
 A. Makhlouf and S. D. Silvestrov,   Hom-algebras and Hom-coalgebras. {\it J. Algebra
Appl.} {\bf 9} (2010) 553-589.


\bibitem{RT93} D. E. Radford and J. Towber, Yetter-Drinfeld categories associated to an arbitrary bialgebra. {\it J. Pure Appl. Algebra}  {\bf 87}(3) (1993) 259-279.

 \bibitem{S69} M. E. Sweedler, \emph{Hopf Algebras},  New York: Benjamin, 1969.


\bibitem{WG} S. X. Wang and S. J. Guo. Symmetries and the u-condition in Hom-Yetter-Drinfeld categories.
{\it J. Math. Phys.} {\bf 55} (2014) 081708.


\bibitem{WCZ12} Z. W. Wang, Y. Y. Chen, L. Y .Zhang. The antipode and Drinfel'd double of Hom-Hopf algebras.
 {\it Scientia Sinica Mathematica.}  {\bf 42}(11) (2012) 1079-1093.

\bibitem{DY09}
D. Yau,  Hom-algebras and homology. {\it J. Lie Theory} {\bf 10} (2009) 409-421.


\bibitem{Y10}
D. Yau, Hom-bialgebras and comodule algebras. {\it Int. Electron. J. Algebra}
{\bf 8} (2010) 45-64.

\bibitem{Y09} D. Yau.  The Hom-Yang-Baxter equation, Hom-Lie algebras, and quasi-triangular bialgebras.
{\em J. Phys. A.} {\bf 42}(16)  (2009) 165202.

\end{thebibliography}
\end{document}